\input{style/arxiv-general.cfg}
\documentclass[MSNbibl,number,citesort,seceqn,dvips]{arxbj}
\makeatletter
   \@ifpackageloaded{graphicx}{}{\usepackage{graphicx}}
\makeatother
\usepackage{upgreek}
\usepackage{graphicx}

\aid{0}
\volume{21}
\issue{2}
\pubyear{2015}
\firstpage{930}
\lastpage{956}
\doi{10.3150/13-BEJ593} 

\makeatletter

\newcommand{\rright}{\right}
\newcommand{\lleft}{\left}
\newcommand{\rrvert}{\vert}
\newcommand{\binom[2]}{\pmatrix{#1\cr #2}}
\renewcommand{\backslash}{\setminus}
\renewcommand{\pi}{\uppi}
\renewcommand{\emptyset}{\varnothing}

\newtheorem{lemma}{Lemma}[section]
\newtheorem{theorem}[lemma]{Theorem}
\newremark{remark}{Remark}[section]
\newremark{Example}{Example}[section]
\newtheorem{Conclusion}{Conclusion}
\newcommand{\ep}{\varepsilon}
\newcommand{\dimp}{\operatorname{dim}_{\mathrm{P}}}
\renewcommand{\dim}{\operatorname{dim}_{\mathrm{H}}}
\newcommand{\R}{\mathrm{R}}
\newcommand{\N}{\mathrm{N}}
\newcommand{\E}{\mathrm{E}}
\renewcommand{\P}{\mathrm{P}}
\makeatother

\begin{document}
\begin{frontmatter}

\title{Exact moduli of continuity for operator-scaling Gaussian random fields}
\runtitle{Exact moduli of continuity for operator-scaling Gaussian fields}

\begin{aug}
\author[A]{\inits{Y.}\fnms{Yuqiang} \snm{Li}\corref{}\thanksref{A}\ead[label=e1]{yqli@stat.ecnu.edu.cn}},
\author[B]{\inits{W.}\fnms{Wensheng} \snm{Wang}\thanksref{B}\ead[label=e2]{wswang@stat.ecnu.edu.cn}} \and
\author[C]{\inits{Y.}\fnms{Yimin} \snm{Xiao}\thanksref{C}\ead[label=e3]{xiao@stt.msu.edu}}
\address[A]{School of Finance and Statistics, East China
Normal University, Shanghai 200241, China.\\ \printead{e1}}
\address[B]{Department of Mathematics,
Hangzhou Normal University, Hangzhou 310036, China.\\ \printead{e2}}
\address[C]{Department of Statistics and Probability,
Michigan State University, 619 Red Cedar Road, East Lansing, MI 48824,
USA. \printead{e3}}
\end{aug}

\received{\smonth{8} \syear{2013}}

%
\begin{abstract}
Let $X=\{X(t), t \in\mathrm{R}^N\}$ be a centered real-valued operator-scaling
Gaussian random field with stationary increments, introduced by Bierm\'{e},
Meerschaert and Scheffler (\textit{Stochastic Process. Appl.}
\textbf{117} (2007) 312--332).
We prove that $X$ satisfies
a form of strong local nondeterminism and establish its exact uniform
and local
moduli of continuity. The main results are expressed in terms of the
quasi-metric $\tau_{E}$ associated with the scaling exponent of $X$.
Examples are provided to illustrate the subtle changes of the regularity
properties.
\end{abstract}

%

\end{frontmatter}

\section{Introduction}

For random fields, ``anisotropy''
is a distinct property from those of one-parameter
processes and is not only important in probability (e.g., stochastic
partial differential
equations) and statistics (e.g., spatio-temporal modeling), but also
for many applied
areas such as economic, ecological, geophysical and medical sciences.
See, for example,
Benson \textit{et al.} \cite{BMB06}, Bonami and Estrade \cite
{BE03}, Chil\'es and
Delfiner \cite{CD99}, Davies and Hall \cite{DH99},
Stein \cite{Stein05,Stein12}, Wackernagel \cite{Wack98}, Zhang \cite
{Zhang07}, and their
combined references for further information.

Many anisotropic random fields $Y = \{Y(t), t \in\R^N\}$ in the
literature have the following scaling property: There exists a linear operator
$E$ (which may not be unique) on $\R^N$ such that for all constants $c
> 0$,
%
%
\begin{equation}
\label{def:OSS} \bigl\{ Y\bigl(c^E t\bigr), t \in\R^N \bigr
\} \stackrel{\mathrm{f.d.}} {=} \bigl\{c Y(t), t \in\R^N \bigr\}.
\end{equation}
Here and in the sequel, ``$ \stackrel{\mathrm{f.d.}}{=} $'' means
equality in all
finite-dimensional distributions and, for $c > 0$, $c^E$ is the linear
operator on $\R^N$ defined by $c^E = \sum_{n=0}^\infty\frac{(\ln
c)^n E^n} {n!}$.
A random field $Y = \{Y(t), t \in\R^N\}$ which satisfies (\ref{def:OSS}),
is called operator-scaling in the time variable (or simply
operating-scaling) with
exponent $E$. Two important examples of real-valued operator-scaling Gaussian
random fields are fractional Brownian sheets introduced by Kamont \cite{K96}
and those with stationary increments introduced by Bierm\'{e},
Meerschaert and
Scheffler \cite{BMS07}. Multivariate random fields with
operator-scaling properties
in both time and space variables have been constructed by Li and Xiao
\cite{LX09}.

Several authors have studied probabilistic and geometric properties of
anisotropic Gaussian random fields. For example, Dunker \cite
{Dunker2000}, Mason
and Shi \cite{MS2001}, Belinski and Linde \cite{BL02}, K\"uhn and
Linde \cite{KL2002}
studied the small ball probabilities of a fractional Brownian sheet
$B^H$, where $H=(H_1,
\ldots, H_N)\in(0, 1)^N$. Mason and Shi \cite{MS2001} also computed
the Hausdorff
dimension of some exceptional sets related to the oscillation of the
sample paths
of $B^H$. Ayache and Xiao \cite{AX05}, Ayache \textit{et al.} \cite
{AWX08}, Wang \cite{Wang07},
Wu and Xiao \cite{WX07}, Xiao and Zhang \cite{XZ04} studied uniform
modulus of
continuity, law of iterated logarithm, fractal properties and joint
continuity of the
local times of fractional Brownian sheets. Wu and Xiao \cite{WX11}
proved sharp
uniform and local moduli of continuity for the
local times of Gaussian fields which satisfy sectorial local nondeterminism.
Luan and Xiao \cite{LuanXiao10} determined the exact Hausdorff measure
functions
for the ranges of Gaussian fields which satisfy strong local nondeterminism.
Meerschaert \textit{et al.} \cite{MWX10} established
exact modulus of continuity for Gaussian fields which satisfy the condition
of sectorial local nondeterminism. Their results and methods are
applicable to
fractional Brownian sheets and certain operator-scaling Gaussian random fields
with stationary increments whose scaling exponent is a diagonal matrix.
We remark that there are subtle differences between certain sample path
properties
of fractional Brownian sheets and those of anisotropic Gaussian random
fields with
stationary increments. This is due to their different properties of
local nondeterminism;
see Xiao \cite{Xiao09} and Li and Xiao \cite{LX13}  for more information.

For an operator-scaling Gaussian random field $X= \{X(t), t \in\R^N\}
$ with stationary increments,
Bierm\'{e} \textit{et al.} \cite{BMS07} showed that the critical
global or directional H\"{o}lder exponents
are given by the real parts of the eigenvalues of the exponent matrix
$E$. The main objective
of this paper is to improve their results and to establish exact
uniform and local moduli of
continuity for these Gaussian fields. Our approach is an extension of
the method in Meerschaert \textit{et al.} \cite{MWX10}.
In particular, we prove in Theorem~\ref{Thm:SLND}
that $X$ has
the property of strong local nondeterminism, which is expressed in
terms of the natural
quasi-metric $\tau_{E}(t-s)$ associated with the scaling exponent $E$ (see
Section~\ref{sec2} for its definition and properties). As an
application of
Theorem~\ref{Thm:SLND} and the
method in \cite{MWX10}, we establish exact uniform and local moduli of
continuity for
$X$ (see Theorems~\ref{Th:UMC} and~\ref{Th:LIL2} below).\looseness=-1

It should be mentioned that Bierm\'{e} \textit{et al.} \cite{BMS07}
constructed
a large class of operator-scaling \mbox{$\alpha$-stable} random fields for
any $\alpha\in(0, 2]$.
By using a LePage-type series representation for stable random fields,
Bierm\'{e} and Lacaux
\cite{BL07} studied uniform modulus of continuity of these
operator-scaling stable random fields.
See also Xiao \cite{Xiao10} for related results using a different
approach based
on the chaining argument. In this paper we will focus on the Gaussian
case (i.e., $\alpha= 2$)
and our Theorem~\ref{Th:UMC} establishes the exact uniform modulus of
continuity, which is more
precise than the results in \cite{BL07} and \cite{Xiao10}.

The rest of this paper is divided into five sections. In Section~\ref{sec2},
we prove some basic properties on the quasi-metric $\tau_{E}$
associated with the scaling exponent $E$ and recall from \cite{BMS07}
the definition of an operator-scaling Gaussian field $X=\{X(t), t \in
\R^N\}$
with stationary increments. In Section~\ref{Sec:SLND}, we prove
the strong local nondeterminism of $X$, and in Sections~\ref{sec4}
and~\ref{sec5}
we
prove the exact uniform and local moduli of continuity of $X$, respectively.
In Section~\ref{sec6} we provide two examples to illustrate our main theorems.

We end the Introduction with some notation. The parameter space is
$\R^N$, endowed with the Euclidean norm
$\|\cdot\|$. For any given two points $s=(s_1,\ldots,s_N)$,
$t=(t_1,\ldots,t_N)$, the inner product of $s,t\in\R^N$ is denoted
by $\langle s, t\rangle$. For $x\in\R_+$, let $\log x:=\ln
(x\vee e
)$ and $\log\log x:=\ln ((\ln x)\vee e )$.
Throughout this paper, we will use $C$ to denote an unspecified positive
and finite constant which may be different in each occurrence. More specific
constants are numbered as $C_1, C_{2},\ldots$\,.

\section{Preliminaries}\label{sec2}

In this section, we show some basic properties of a real $N\times N$
matrix $E$ and prove several lemmas on
the quasi-metric $\tau_{E}$ on $\R^N$. Then we recall from Bierm\'
{e} \textit{et al.} \cite{BMS07} the definition of operator-scaling Gaussian
random fields with a harmonizable representation.

For a real $N\times N$ matrix $E$, it is well known that
$E$ is similar to a real Jordan canonical form, i.e. there exists a
real invertible $N\times N$ matrix $P$ such that
\[
E=PDP^{-1},
\]
where $D$ is a real $N\times N$ matrix of the form
%
%
\begin{equation}
\label{Eq:J0} D= \pmatrix{ 
J_1& 0&
\cdots&0
\cr
0 &J_2& \cdots&0
\cr
\vdots&\vdots&\ddots&\vdots
\cr
0&
0& \cdots&J_p 
}
\end{equation}
and $J_i$, $1\leq i\leq p$, is either Jordan cell matrix of the form
%
%
\begin{equation}
\label{Eq:J1} \pmatrix{ \lambda&0& 0& \cdots&0
\cr
1& \lambda&0& \cdots&0
\cr
0&
1& \lambda&\cdots&0
\cr
\vdots&\vdots&\vdots&\ddots&\vdots
\cr
0& 0& 0& \cdots&
\lambda }
\end{equation}
with $\lambda$ a real eigenvalue of $E$ or blocks of the form
%
%
\begin{equation}
\label{Eq:J2} \lleft( %
\begin{array} {@{}c@{\quad}c@{\quad}c@{\quad}c@{
\quad}c@{}} \Lambda&0& 0 &\cdots&0
\\
I_2& \Lambda&0& \cdots&0
\\
0& I_2& \Lambda&\cdots&0
\\
\vdots& \vdots& \vdots& \ddots& \vdots
\\
0 &0& 0& \cdots&\Lambda \end{array} %
\rright) \qquad\mbox{with }
I_2=\lleft( %
\begin{array} {@{}c@{\quad}c@{}} 1 &0
\\
0 &1 \end{array} %
\rright) \mbox{ and } \Lambda=\lleft( %
\begin{array} {@{}c@{\quad}c@{}} a &-b
\\
b &a \end{array} %
\rright),
\end{equation}
where the complex numbers $a\pm ib$ ($b\neq0$) are complex conjugated
eigenvalues of $E$.

Denote the size of $J_k$ by $\tilde{l}_k$ and let $a_k$ be the real
part of the corresponding eigenvalue(s) of $J_k$. Throughout this
paper, we always suppose that
\[
1<a_1\leq a_2\leq\cdots\leq a_p.
\]
Note that $p\leq N$,
$\tilde{l}_1+\tilde{l}_2+\cdots+\tilde{l}_p=N$ and
$Q:=\operatorname{trace}(E)=\sum_{j=1}^p a_j\tilde{l}_j$.

As done in Bierm\'e and Lacaux \cite{BL07}, we can construct the
$E$-invariant subspace $W_k$ associated with $J_k$ by
\[
W_k=\operatorname{span} \Biggl\{f_j \dvt \sum
_{i=1}^{k-1}\tilde {l}_i+1\leq j\leq
\sum_{i=1}^k\tilde{l}_i
\Biggr\},
\]
where $f_j$ is the $j$th column vector of the matrix $P$. Then $\R^N$
has a direct sum decomposition of
\[
\R^N = W_1 \oplus\cdots\oplus W_p.
\]

It follows from Meerschaert and Scheffler \cite{RVbook}, Chapter~6
(see also \cite{BMS07}, Section~2)
that there exists a
norm $\|\cdot\|_E$ on $\R^N$ such that for the unit sphere $S_E =
\{x \in\R^N\dvt \|x\|_E=1 \}$
the mapping $\Psi\dvtx (0,\infty) \times S_E \to\R^{N}\backslash\{
0\}$
defined by $\Psi(r,\theta) = r^E \theta$
is a homeomorphism. Hence, every $x\in\R^N\setminus\{0\}$ can be
written uniquely as $x=(\tau_{E}(x))^E
l_E(x)$ for some radial part $\tau_{E}(x)>0$ and some direction
$l_E(x)\in S_E$ such that the functions $x\mapsto
\tau_{E}(x)$ and $x\mapsto l_E(x)$ are continuous. For $x \in\R
^N\backslash\{0\}$,
$ (\tau_{E}(x), l_E(x) )$ is referred to as its polar
coordinates associated with $E$.

It is shown in \cite{RVbook} that $\tau_{E}(x)=\tau_{E}(-x)$ and
$\tau_{E}(r^E x)=r\tau_{E}(x)$ for all
$r>0$ and $x\in\R^N\backslash\{0\}$. Moreover, $\tau_{E}(x)\to
\infty$ as
$x\to\infty$ and $\tau_{E}(x)\to0$ as $x\to0$. Hence, we can
extend $\tau_{E}(x)$
continuously to $\R^N$ by setting $\tau_{E}(0) = 0$.

The function $\tau_{E}(x)$ will play essential roles in this paper.
We first recall some known
facts about it.
\begin{enumerate}[(iii)]
\item[(i)] Lemma~2.2 in \cite{BMS07} shows that there
exists a constant $C \geq1$ such that
%
%
\begin{equation}
\label{s2-6} \tau_{E}(x+y)\leq C \bigl( \tau_{E}(x) +
\tau_{E}(y) \bigr),\qquad \forall x, y \in \R^N.
\end{equation}
Hence, we can regard $ \tau_{E}(x-y)$ as a quasi-metric on $\R^N$.

\item[(ii)] Since the norms $\|\cdot\|_E$ and $\|\cdot\|$ are
equivalent, Lemma~2.1 in
\cite{BMS07} implies that for any $0 < \delta< a_1$ there exist finite
constants $C_{1}, C_{2} > 0$, which may depend on $\delta$, such that
for all
$\|x\|\le1$ or all $\tau_{E}(x)\le1$,
%
%
\begin{equation}
\label{s2-1} C_{1}\|x\|^{1/(a_1-\delta)}\leq\tau_{E}(x)\leq
C_{2}\|x\| ^{1/(a_p+\delta)},
\end{equation}
and, for all $\|x\|>1$ or all $\tau_{E}(x)>1$,
%
%
\begin{equation}
\label{s2-2} C_{1}\|x\|^{1/(a_p+\delta)}\leq\tau_{E}(x)\leq
C_{2}\|x\| ^{1/(a_1-\delta)}.
\end{equation}

\item[(iii)] Bierm\'e and Lacaux \cite{BL07}, Corollary~3.4, proved
the following improvement of
(\ref{s2-1}): For any $\eta\in(0,1)$, there exists a finite constant
$C_{3} \ge1$ such that for all $x\in W_j\setminus\{0\}$,
$1\leq j\leq p$, with $\|x\|\leq\eta$
%
%
\begin{equation}
\label{s2-3} C_{3}^{-1}\|x\|^{1/a_j} \bigl |\ln\|x\|
\bigr |^{-(l_j-1)/a_j}\leq\tau _{E}(x)\leq C_{3}\|x
\|^{1/a_j} \bigl |\ln\|x\| \bigr |^{(l_j-1)/a_j},
\end{equation}
where $l_k=\tilde{l}_k$ if $J_k$ is a Jordan cell matrix as in (\ref
{Eq:J1}) or $l_k=\tilde{l}_k/2$
if $J_k$ is of the form (\ref{Eq:J2}). 
\end{enumerate}

We remark that, as shown by Example~\ref{exa6.2} below, both the upper
and lower
bounds in (\ref{s2-3})
can be achieved and this fact makes the regularity properties of an
operator-scaling Gaussian field
with a general exponent $E$ more intriguing.

For any $x \in\R^N$, let $x = \bar{x}_1 \oplus\bar{x}_2\oplus
\cdots\oplus\bar{x}_p$
be the direct sum decomposition of $x$ in the $E$-invariant subspaces
$W_j$, $j=1,2,\ldots,p$.
This notation is used in Lemmas~\ref{lem-s2-3} and~\ref{lem-s2b}.

%
\begin{lemma}\label{lem-s2-3}
There exists a finite constant $C>0$ such that for all $x \in\R^N$
and $j=1,2,\ldots, p$,
we have
%
%
\begin{equation}
\label{Eq:tau1} \tau_{E}(\bar{x}_j)\leq C
\tau_{E}(x).
\end{equation}
\end{lemma}

\begin{pf}
Since (\ref{Eq:tau1}) holds trivially for $x = 0$. We
only consider
$x\in\R^N\backslash\{0\}$, which can be written as $x=(\tau_E(x))^E l_E(x)$
for some $l_E(x)\in S_E$. Denote the direct sum decomposition of
$l_E(x)$ in the $E$-invariant
subspaces $W_j$, $j=1,2,\ldots,p$, by $l_E(x) = x_1'\oplus\cdots
\oplus x_p'$. Then
from the fact that $ (\tau_E(x) )^Ex_j'\in W_j$ for all
$j=1,2,\ldots,p$, it follows that
\[
\bar{x}_j=\bigl(\tau_E(x)\bigr)^E
x_j'.
\]
Since $S_E$ is bounded, that is, there exists $M>0$ such that
$S_E\subset\{y \in\R^N\dvt\|y\|\leq M\}$,
we can easily see that $x_j'\in\{y \in\R^N\dvt\|y\|\leq M\}$ for all
$j=1,2,\ldots,p$.
Let $C=\max_{\|x\|\leq M}\tau_E(x)\in(0, \infty)$. Then for all
$j=1,2,\ldots, p$
\[
\tau_E(\bar{x}_j)=\tau_E(x)
\tau_E\bigl(x_j'\bigr)\leq C
\tau_E(x),
\]
which is the desired conclusion.
\end{pf}

As a consequence of (\ref{s2-6}) and Lemma~\ref{lem-s2-3}, we have
the following lemma.

%
\begin{lemma}\label{lem-s2b}
There is a finite constant $C \ge1$ such that
%
%
\begin{equation}
\label{Eq:tau1b} C^{-1}\sum_{i=1}^p
\tau_{E}(\bar{x}_i)\leq\tau_{E}(x) \le C \sum
_{i=1}^p \tau_{E}(
\bar{x}_i), \qquad\forall x \in\R^N.
\end{equation}
\end{lemma}

The following lemma implies that the function $\tau_{E}(x)$ is
O-regular varying
at both the origin and the infinity (cf. Bingham \textit{et al.}
\cite{BGT87}, pages 65--67).

%
\begin{lemma}\label{lem-s2-2}
Give any constants $0 < a < b< \infty$, there exists a finite constant
$C_{4} \ge1$ such that for all
$x\in\R^N$ and $\beta\in[a, b]$,
%
%
\begin{equation}
\label{Eq:tau2} C_{4}^{-1}\tau_{E}(x)\leq
\tau_{E}(\beta x)\leq C_{4}\tau_{E}(x).
\end{equation}
\end{lemma}

\begin{pf} To prove the left inequality in (\ref{Eq:tau2}), note that
$\Lambda=\{\beta x\dvt x\in S_E, \beta\in[a, b]\}$ is a compact set which
does not contain $0$. This and the continuity of $\tau_{E}(\cdot)$
on $\R^N$,
imply $\min_{x\in\Lambda} \tau_{E}(x)>0$. Hence, by taking
$C_{4}^{-1} = 1 \wedge\min_{x\in\Lambda} \tau_{E}(x)$, we have
\[
\tau_{E}(\beta x)=\tau_{E}\bigl(\beta\tau_{E}^E(x)l_E(x)
\bigr) =\tau_{E}(x)\tau_{E}\bigl(\beta l_E(x)
\bigr)\geq C_{4}^{-1} \tau_{E}(x).
\]
The right inequality in (\ref{Eq:tau2}) can be proved in the same way. This
finishes the proof.
\end{pf}

%
\begin{lemma}\label{lem-s2-1}
There is a subsequence $\{n_k\}_{k\in\N} \subseteq\N$ such
that $n_k\geq k$ for all $k \ge1$ and
%
%
\begin{equation}
\label{thm-5} \min_{1\leq i\leq
2^{n_k}}\tau_{E}\bigl(\bigl\langle
i2^{-n_k}\bigr\rangle\bigr)\geq C_{4}^{-1}\tau
_{E}\bigl(\bigl\langle 2^{-n_k}\bigr\rangle\bigr),
\end{equation}
where $\langle c\rangle=(c,c,\ldots,c)\in\R^N$ for any $c\in\R$.
\end{lemma}

\begin{pf}
Suppose $\min_{1\leq i\leq2^{n}}\tau_{E}(\langle i2^{-n}\rangle)$ is
attained at $i=K_n$.
There is an integer $m_n \in[0, n] $ such that $2^{n-m_n-1}<K_n\leq
2^{n-m_n}$.
Therefore, we can rewrite $K_n2^{-n}$ as $\beta2^{-m_n}$ for some
$\beta\in(1/2,1]$. Since
$
\{i2^{-m_n}, i=1,\ldots, 2^{m_n}\}\subset\{i2^{-n}, i=1,\ldots,
2^{n}\}$,
we have
\[
\min_{1\leq i\leq2^{m_n}}\tau_{E}\bigl(\bigl\langle
i2^{-m_n}\bigr\rangle\bigr)\geq \min_{1\leq i\leq
2^{n}}
\tau_{E}\bigl(\bigl\langle i2^{-n}\bigr\rangle\bigr) =
\tau_{E}\bigl(\bigl\langle\beta 2^{-m_n}\bigr\rangle\bigr) \geq
C_{4}^{-1}\tau_{E}\bigl(\bigl
\langle2^{-m_n}\bigr\rangle\bigr),
\]
where the last inequality follows from Lemma~\ref{lem-s2-2} with $[a,
b] = [1/2, 1]$.
Furthermore, by the fact
\[
\min_{1\leq i\leq
2^{n}}\tau_{E}\bigl(\bigl\langle
i2^{-n}\bigr\rangle\bigr)\leq\tau_{E}\bigl(\bigl\langle
2^{-n}\bigr\rangle \bigr)\to0,
\]
as $n\to\infty$, we know that $\tau_{E}(\langle2^{-m_n}\rangle)\to0$
which implies that
$m_n\to\infty$ as $n\to\infty$. Hence, a desired subsequence $\{
n_k\}_{k\in\N}$ can be
selected from $\{m_n\}$.
\end{pf}

Let $E'$ be the transpose of $E$. An $E'$-homogeneous function $\psi
\dvtx
\R^N \to[0,\infty)$ is a function which satisfies that $\psi(x)> 0$
and $\psi(r^{E'}x)=r\psi(x)$ for all $r
> 0$ and $x \in\R^N \backslash\{0\}$. For any continuous
$E'$-homogeneous function
$\psi\dvtx \R^N \to[0,\infty)$, Bierm\'{e} \textit{et al.} \cite
{BMS07}, Theorem~4.1, showed that the real-valued Gaussian random field $X_\psi=
\{X_\psi(t), t \in\R^N\}$, where
%
%
\begin{equation}
\label{Hfield} X_\psi(t)= \operatorname{Re} \int_{\R^N}
\bigl(\mathrm {e}^{\mathrm{i}\langle t,\xi\rangle
}-1 \bigr) \frac{\widetilde{\mathcal{M}} (\mathrm{d}\xi)} { \psi(\xi)^{1 +
Q/2}}, \qquad t\in
\R^N,
\end{equation}
is well defined and stochastic continuous if and only if $\min_{1\le
j \le p}a_j > 1$. In the latter case, they further proved that
$X_\psi$ satisfies (\ref{def:OSS}) and has stationary increments.
Here, $\widetilde{\mathcal M}$ is a centered complex-valued Gaussian
random measure in $\R^N$ with the Lebesgue measure $m_N$ as its
control measure. Namely, $\widetilde{\mathcal{M}}$ is a centered
complex-valued Gaussian process defined on the family $\mathcal{A} =
\{A \subset\R^N\dvt m_N(A) < \infty\}$ which satisfies
%
%
\begin{equation}
\label{Eq:M} \E \bigl(\widetilde{\mathcal M}(A) \overline{ \widetilde{\mathcal
M} (B)} \bigr) = m_N (A \cap B) \quad\mbox{and}\quad \widetilde{
\mathcal M}(-A) = \overline{\widetilde{\mathcal M}(A)}
\end{equation}
for all $A, B \in\mathcal{A}$.

\begin{remark}\label{rem2.1} The following are some remarks on the
Gaussian random field $X_\psi$.
\begin{itemize}
\item If, in addition, $\psi$ is symmetric in the sense that $\psi
(\xi) = \psi(-\xi)$ for all $\xi\in\R^N$,
then because of (\ref{Eq:M}) the Wiener-type integral in the
right-hand side of (\ref{Hfield}) is
real-valued. Thus, in this latter case, ``$\operatorname{Re}$'' in
(\ref
{Hfield}) is not needed. For simplicity, we
assume that $\psi$ is symmetric in the rest of the paper. A large
class of continuous, symmetric
$E'$-homogeneous functions has been constructed in
\cite{BMS07}, Theorem~2.1.

\item By replacing $\widetilde{\mathcal M}$ in (\ref{Hfield}) by a
complex-valued isotropic
$\alpha$-stable random measure $\widetilde{\mathcal M}_\alpha$ with
Lebesgue control measure
(see \cite{ST94}, page 281), Bierm\'{e} \textit{et al.} \cite
{BMS07}, Theorem~3.1, obtained a class of
\emph{harmonizable} operator-scaling $\alpha$-stable random fields.
They also defined a
class of operator-scaling $\alpha$-stable fields by using
moving-average representations. When
$\alpha\in(0, 2)$, stable random fields with harmonizable and
moving-average representations
are generally different. However, for the Gaussian case of $\alpha=2$,
the Planchrel theorem
implies that every Gaussian random field with a moving-average
representation in \cite{BMS07}
also has a harmonizable representation of the form (\ref{Hfield}).
\end{itemize}
\end{remark}

\section{Strong local nondeterminism of operator-scaling Gaussian~fields}\label{Sec:SLND}

Let $E$ be an $N\times N$ matrix such that the real parts of its eigenvalues
satisfy $\min_{1\le j \le p}a_j > 1$ and let $\psi$ be a continuous,
symmetric,
$E'$-homogeneous function with $\psi(x)> 0$ for $x \neq0$ as in
Section~\ref{sec2}.
Let $X_{\psi} = \{X_\psi(t), t \in\R^N\}$ be the operator-scaling
Gaussian field
with scaling exponent $E$, defined by (\ref{Hfield}). For simplicity, we
write $X_\psi$ as $X$. Note that the assumptions on $\psi$ imply
%
%
\begin{equation}
\label{Eq:psibound} 0<m_{\psi}=\min_{x\in S_{E'}}\psi(x)\leq
M_{\psi} =\max_{x\in S_{E'}}\psi(x)<\infty.
\end{equation}

The dependence structure of the operator-scaling Gaussian field
$X$ is complicated for a general matrix $E$. In order to
study sample path properties and characterize the anisotropic
nature of $X$, we prove that $X$ has the property
of ``strong local nondeterminism'' with respect to the quasi-metric
$\tau_{E}(s-t)$. The main result of this section is Theorem~\ref{Thm:SLND}, which extends Theorem~3.2 in Xiao \cite{Xiao09} and
will play an important role in Section~\ref{sec4} below.

Since many sample path properties of $X$ are determined by the
canonical metric\vspace*{1pt}
%
%
\begin{equation}
\label{s3-0} d_X(s,t)= \bigl[\E \bigl(X(s)-X(t)
\bigr)^2 \bigr]^{1/2},\qquad\forall s,t\in\R^N,
\end{equation}
our first step is to establish the relations between $d_X(s, t)$ and
$\tau_{E}(s-t)$.

%
\begin{lemma}\label{Lem:Var1}
There exists a finite constant $C \ge1$ such that\vspace*{1pt}
%
%
\begin{equation}
\label{Eq:Var1} C^{-1} \tau_{E}^2(s-t)\leq
d_{X}^2(s, t)\leq C \tau_{E}^2(s-t),
\qquad\forall s, t\in\R^N.
\end{equation}
\end{lemma}

\begin{pf}
For all $s, t\in\R^N$, by (\ref{Hfield}), we have
\begin{eqnarray*}
d_{X}^2(s,t)&=&\int_{\R^N} \bigl |
\mathrm{e}^{\mathrm{i}\langle
s,x\rangle
}-\mathrm{e}^{\mathrm{i}\langle t, x\rangle} \bigr |^2
\frac{\mathrm{d} x}{\psi(x)^{2+Q}}
\\
&=&2\int_{\R^N}\bigl(1-\cos\langle s-t,x\rangle\bigr)
\frac{\mathrm{d}
x}{\psi(x)^{2+Q}}.
\end{eqnarray*}
Let $y=\tau_{E}^{E'}(s-t)x$. Then $\mathrm{d}x=(1/\tau_{E}(s-t))^Q \,
\mathrm{d}y$. Hence,
%
%
\begin{equation}
\label{s3-1} d_{X}^2(s,t)=2\tau_{E}(s-t)^2
\int_{\R^N} \biggl(1-\cos \biggl\langle \biggl(
\frac{1}{\tau_{E}(s-t)} \biggr)^E(s-t),y \biggr\rangle \biggr)
\frac{\mathrm{d}
y}{\psi(y)^{2+Q}}.
\end{equation}
Since for all $s \ne t$, $\tau_{E} ( (\frac{1}{\tau
_{E}(s-t)} )^E(s-t) )=1$.
Hence, the set
\[
\biggl\{ \biggl(\frac{1}{\tau_{E}(s-t)} \biggr)^E(s-t)\dvt s \ne t \in \R
^N \biggr\}
\]
is compact and does not contain $0$. On the other hand, a slight
modification of
the proof of Theorem~4.1 in \cite{BMS07} shows that the
function $\xi\mapsto\int_{\R^N} ( 1 - \cos\langle\xi,y\rangle
)\frac{\mathrm{d}
y}{\psi(y)^{2+Q}}$ is continuous on $\R^N$ and positive on $\R
^N\backslash\{0\}$.
Therefore, the last
integral in (\ref{s3-1}) is bounded from below and above by positive
and finite
constants. This proves (\ref{Eq:Var1}).
\end{pf}

%
\begin{theorem}\label{Thm:SLND}
There exists a constant $C_5> 0$ such that for all $n \geq2$ and all $
t^1, \ldots, t^n \in\R^N$, we have
%
%
\begin{equation}
\label{eq:slnd2} \operatorname{Var} \bigl( X\bigl(t^n\bigr) | X
\bigl(t^1\bigr), \ldots, X\bigl(t^{n-1}\bigr) \bigr) \ge
C_5 \min_{0 \le k \le n-1} \tau _{E}^2
\bigl(t^n-t^k\bigr),
\end{equation}
where $t^0 =0$.
\end{theorem}

\begin{pf} The proof is a modification of that of Theorem~3.2 in
Xiao \cite{Xiao09}. We denote $r=\min_{0\leq k\leq n-1}\tau
_{E}(t^n-t^k)$. Since
\[
\operatorname{Var} \bigl(X\bigl(t^n\bigr) |X(t_1),
\ldots,X\bigl(t^{n-1}\bigr) \bigr) =\inf_{u_1,\ldots,u_{n-1}\in\R} \E
\Biggl[ \Biggl(X\bigl(t^n\bigr)-\sum_{k=1}^{n-1}u_kX
\bigl(t^k\bigr) \Biggr)^2 \Biggr],
\]
it suffices to prove the existence of a constant $C > 0$ such that
\begin{equation}
\label{s3-2}
\E \Biggl[ \Biggl(X\bigl(t^n\bigr)-\sum
_{k=1}^{n-1}u_kX\bigl(t^k
\bigr) \Biggr)^2 \Biggr]\geq C r^2
\end{equation}
for all $u_k\in\R$, $k=1,2,\ldots,n-1$. It follows from (\ref
{Hfield}) that
\[
\E \Biggl[ \Biggl(X\bigl(t^n\bigr)-\sum
_{k=1}^{n-1}u_kX\bigl(t^k
\bigr) \Biggr)^2 \Biggr]= \int_{\R^N} \Biggl|
\mathrm{e}^{\mathrm{i}\langle t^n, x\rangle}-\sum_{k=0}^{n-1}u_k
\mathrm{e}^{\mathrm{i}\langle
t^k,x\rangle} \Biggr|^2 \frac{\mathrm{d}x}{\psi(x)^{2+Q}},
\]
where $t^0=0$ and $u_0=1-\sum_{k=1}^n u_k$. Let
$\delta(\cdot)\dvtx\R^N\mapsto[0,1]$ be a function in $C^\infty
(\R^N)$
such that $\delta(0)=1$ and it vanishes outside the open set
$B=\{x\dvt\tau_{E}(x)<1\}$. Denote by $\widehat{\delta}$ the Fourier
transform of $\delta$. Then $\widehat{\delta}\in C^\infty(\R^N)$ as
well and decays rapidly as $x\to\infty$, that it, for all $\ell\ge
1$, we have
$\|x\|^\ell|\widehat{\delta}(x)|\to0$ as $x\to\infty$. This and
(\ref{s2-2})
further imply that for all $\ell\ge1$,
%
%
\begin{equation}
\label{s3-3} \tau_{E}(x)^\ell\bigl |\widehat{\delta}(x)\bigr |\to0
\qquad\mbox{as } x\to \infty.
\end{equation}
Let $\delta_r(t)=r^{-Q}\delta(r^{-E}t)$. Then
\[
\delta_r(t)=(2\pi)^{-N}\int_{\R^N}
\mathrm{e}^{-\mathrm{i}\langle
t, x\rangle
}\widehat{\delta}\bigl(r^{E'}x\bigr)\,\mathrm{d}
x.
\]
Since $\min\{\tau_{E}(t^n-t^k),0\leq k\leq n-1\}=r$, we have
$\delta_r(t^n-t^k)=0$ for all $k=0,1,\ldots,n-1$. Therefore,
%
%
\begin{eqnarray}
\label{s3-4} J&:=&\int_{\R^N} \Biggl(\mathrm{e}^{\mathrm{i}\langle t^n,x\rangle
}-
\sum_{k=0}^{n-1}\mathrm{e}^{\mathrm{i}\langle
t^k,x\rangle}
\Biggr)\mathrm{e}^{-\mathrm{i}\langle t^n,x\rangle
}\widehat{\delta }\bigl(r^{E'}x\bigr)\,
\mathrm{d}x
\nonumber
\\[-8pt]
\\[-8pt]
&=&(2\pi)^N \Biggl(\delta_r(0)-\sum
_{k=0}^{n-1}u_k\delta _r
\bigl(t^n-t^k\bigr) \Biggr) =(2\pi)^N
r^{-Q}.
\nonumber
\end{eqnarray}
By H\"{o}lder's inequality, a change of variables,
the $E'$-homogeneity of $\psi$ and (\ref{Eq:psibound}), we derive
%
%
\begin{eqnarray}
\label{s3-5} J^2&\leq&\int_{\R^N}\Biggl |
\mathrm{e}^{\mathrm{i}\langle
t^n,x\rangle}-\sum_{k=0}^{n-1}
\mathrm{e}^{\mathrm{i}\langle
t^k,x\rangle} \Biggr|^2\frac{\mathrm{d}x}{\psi(x)^{2+Q}}\int
_{\R
^N}\psi(x)^{2+Q} \bigl |\widehat{\delta}
\bigl(r^{E'}x\bigr)\bigr |^2\,\mathrm{d}x
\nonumber
\\
&=& r^{-2Q-2}\E \Biggl(\Biggl |X\bigl(t^n\bigr)-\sum
_{k=1}^{n-1}u_k X\bigl(t^k
\bigr)\Biggr |^2 \Biggr) \int_{\R^N}\psi(y)^{2+Q}\bigl |
\widehat{\delta}(y)\bigr |^2\,\mathrm {d}y
\nonumber
\nonumber
\\[-8pt]
\\[-8pt]
&\leq&r^{-2Q-2}\E \Biggl( \Biggl|X \bigl(t^n\bigr)-\sum
_{k=1}^{n-1}u_k X\bigl(t^k
\bigr) \Biggr|^2 \Biggr) \int_{\R^N}\tau_{E}(y)^{2+Q}M_\psi^{2+Q}
\bigl |\widehat{\delta }(y) \bigr |^2 \,\mathrm{d}y
\nonumber
\\
&\leq& Cr^{-2Q-2}\E \Biggl( \Biggl|X\bigl(t^n\bigr)-\sum
_{k=1}^{n-1}u_k X\bigl(t^k
\bigr) \Biggr|^2 \Biggr)
\nonumber
\end{eqnarray}
for some finite constant $C>0$, since
$\int_{\R^N}\tau_{E}(y)^{2+Q}|\widehat{\delta}(y)|^2\,\mathrm
{d}y<\infty$
which follows from (\ref{s3-3}). Combining (\ref{s3-4}) and (\ref{s3-5})
yields
(\ref{s3-2})
for an appropriate constant $C_5>0$.
\end{pf}

The relation (\ref{eq:slnd2}) is a property of strong local
nondeterminism, which is more general than that in Xiao \cite{Xiao09}
and can be applied to establish many sample path properties of~$X$.

For any $s, t \in\R^N$ with $s \ne t$, we decompose $s-t$ as a direct
sum of
elements in the $E$-invariant subspaces $W_j$, $j=1,2,\ldots,p$,
\[
s-t = (s_1-t_1)\oplus\cdots\oplus(s_p-t_p).
\]
Then (\ref{s2-3}) and Lemmas~\ref{Lem:Var1} and~\ref{lem-s2b} imply
%
%
\begin{eqnarray}
\label{Eq:Var2} 
&& C^{-1}\sum
_{j=1}^p \|s_j-t_j
\|^{1/a_j} \bigl |\ln\|s_j - t_j\| \bigr |^{-(l_j-1)/a_j}\nonumber
\\
&&\quad 
\leq  d_{X}^2(s, t)
\\
&&\quad 
\leq C \sum_{j=1}^p
\|s_j-t_j\|^{1/a_j} \bigl |\ln\|s_j -
t_j\| \bigr |^{(l_j-1)/a_j}.\nonumber
\end{eqnarray}
Moreover, Theorem~\ref{Thm:SLND} implies that for all $n \geq2$ and
all $
t^1, \ldots, t^n \in\R^N$, we have
%
%
\begin{eqnarray}
\label{eq:slnd3} %
&&\operatorname{Var} \bigl(   X
\bigl(t^n\bigr) \rrvert X\bigl(t^1\bigr), \ldots, X
\bigl(t^{n-1}\bigr) \bigr)
\nonumber
\\[-8pt]
\\[-8pt]
&&\quad\ge C \min_{0 \le k \le n-1} \sum_{j=1}^p
\bigl \|t^n_j-t^k_j
\bigr \|^{1/a_j} \bigl |\ln\bigl \|t^n_j - t^k_j
\bigr \|\bigr  |^{-(l_j-1)/a_j},
\nonumber
\end{eqnarray}
where $t^0 =0$.

Inequalities (\ref{Eq:Var2}) and (\ref{eq:slnd3}) are similar to
Condition (C1) and (C$3'$)
in Xiao \cite{Xiao09}. Hence, many results on the Hausdorff dimensions
of various random sets
and joint continuity of the local times can be readily derived from
those in \cite{Xiao09},
and these results can be explicitly expressed in terms of the real
parts $\{a_j, 1\le j \le p\}$ of
the eigenvalues of the scaling exponent $E$.

To give some examples,
we define a vector $(H_1, \ldots, H_N) \in(0, 1)^N$ as follows.

For $1 \le i \le\tilde{l}_p$, define $H_{i} = a_p^{-1}$. In
general, if
$1 +\sum_{j=k}^p \tilde{l}_j \le i \le\sum_{j=k-1}^p \tilde{l}_j$
for some $ 2 \le k \le p$,
then we define $H_{i} = a_{k-1}^{-1}$. Since $1 < a_1 \le a_2 \le
\cdots\le a_p$, we have
\[
0 < H_1 \le H_2 \le\cdots\le H_N < 1.
\]

Consider a Gaussian random field $\vec X = \{\vec X(t), t \in\R^N\}$
with values in $\R^d$ defined by
\[
\vec X(t) = \bigl(X_1(t), \ldots, X_d(t)\bigr),
\]
where $X_1, \ldots, X_d$ are independent copies of the centered
Gaussian field $X$ in the above.
Let $\vec X([0, 1]^N)$ and $\operatorname{Gr}\vec X([0, 1]^N) = \{(t,
\vec X(t)), t
\in[0, 1]^N\}$ denote respectively the
range and graph of $\vec X$, then Theorem~6.1
in \cite{Xiao09} implies that with probability 1,
%
%
\begin{equation}
\label{Eq:range} \dim{\vec X} \bigl([0, 1]^N \bigr) = \dimp{\vec X}
\bigl([0, 1]^N \bigr) = \min \Biggl\{ d; \sum
_{j=1}^N \frac{1} {H_j} \Biggr\},
\end{equation}
where $\dim$ and $\dimp$ denote Hausdorff and packing dimension
respectively, and
%
%
\begin{eqnarray}
\label{Eq:graph} 
 \dim\operatorname{Gr} {\vec X} \bigl([0,
1]^N \bigr) &=& \dimp \operatorname{Gr} {\vec X} \bigl([0,
1]^N \bigr)
\nonumber
\\
&  =& \min \Biggl\{ \sum_{j=1}^k
\frac{H_k} {H_k} + N-k + (1 - H_k) d, 1 \le k \le N; \sum
_{j=1}^N \frac{1} {H_j} \Biggr\}
\\
& =& \cases{ \displaystyle\sum_{j=1}^N
\frac{1} {H_j}, &\quad$\mbox{if } \displaystyle\sum
_{j=1}^N \frac{1} {H_j} \le d$,
\cr
\displaystyle\sum_{j=1}^k
\frac{H_k} {H_j} + N-k + (1 - H_k) d, &\quad$\mbox{if }
\displaystyle\sum_{j=1}^{k-1}
\frac{1} {H_j} \le d < \sum_{j=1}^{k}
\frac{1} {
H_j} $, }
\nonumber
\end{eqnarray}
where $\sum_{j=1}^{0} \frac{1} {H_j} := 0$.

Similarly, it follows from Theorem~7.1 in Xiao \cite{Xiao09} that the
following hold:
\begin{enumerate}[(ii)]
\item[(i)] If $\sum_{j=1}^N \frac{1} {H_j} < d$, then for every $x\in
\R^d$, $\vec X^{-1}(\{x\}) = \emptyset$ a.s.

\item[(ii)] If $ \sum_{j=1}^N \frac{1} {H_j} > d$, then for every
$x\in\R^d$,
%
%
\begin{eqnarray}
\label{Eq:dimlevel} 
\dim\vec X^{-1}\bigl(\{x\}\bigr)&=&\dimp\vec
X^{-1}\bigl(\{x\}\bigr)
\nonumber
\\
&=& \min \Biggl\{ \sum_{j=1}^k
\frac{H_k} {H_j} + N-k - H_k d, 1 \le k \le N \Biggr\}
\\
&=&\sum_{j=1}^k \frac{H_k} {H_j} + N-k
- H_k d,\qquad\mbox{if } \sum_{j=1}^{k-1}
\frac{1} {H_j} \le d < \sum_{j=1}^{k}
\frac{1} {
H_j}
\nonumber
\end{eqnarray}
holds with positive probability.
\end{enumerate}

In light of the dimension results (\ref{Eq:range})--(\ref
{Eq:dimlevel}), it would be interesting
to determine the exact Hausdorff (and packing) measure functions for
the above random sets.
In the special case of fractional Brownian motion, the corresponding
problems have been
investigated by Talagrand \cite{T95}, Xiao \cite{X97a,X97b}, Baraka
and Mountford \cite{BM11}.
For anisotropic Gaussian random fields, the problems are more
difficult. Only an exact
Hausdorff measure function for the range has been determined by Luan
and Xiao
\cite{LuanXiao10} for a special case of anisotropic Gaussian random fields.

\section{Uniform modulus of continuity}\label{sec4}

In this section, we establish the exact modulus of continuity for
$X$. We first rewrite Lemma~7.1.1 in Marcus and Rosen \cite{MR-book}
as follows.

%
\begin{lemma}\label{Lem-s4-1}
Let $\{G(u), u\in\R^N\}$ be a centered Gaussian random field. Let
$\omega\dvtx
\R_+\to\R_+$ be a function with $\omega(0+) = 0$ and $\Gamma
\subset\R^N$
be a compact set. Assume that there is a continuous map $\tau\dvtx
\R^N \mapsto\R_+$ with $\tau(0)=0$ such that $d_G$ is
continuous on $\tau$, i.e., $\tau(u_n-v_n)\to0$ implies
$d_{G}(u_n,v_n)\to0$. Then
\[
\lim_{\delta\to0}\mathop{\sup_{\tau(u-v)\leq\delta}}_{ u,v\in
\Gamma}
\frac{|G(u)-G(v)|}{\omega(\tau(u-v))}\leq C, \qquad\mbox{a.s. for some constant } C<\infty
\]
implies that
\[
\lim_{\delta\to0}\mathop{\sup_{\tau(u-v)\leq\delta}}_{ u,v\in
\Gamma}
\frac{|G(u)-G(v)|}{\omega(\tau(u-v))}=C', \qquad\mbox{a.s. for some constant }
C'<\infty.
\]
This result is also valid for the local modulus of continuity of $G$,
that is, it
holds with $v$ replaced by $u_0$ and with the supremum taken over $u\in
\Gamma$.
\end{lemma}

%
\begin{remark}
Lemma~\ref{Lem-s4-1} is slightly different from Lemma~7.1.1 in Marcus and
Rosen \cite{MR-book}, where $\tau$ is assumed to be a pseudo-norm. However,
by carefully checking its proof in \cite{MR-book}, this requirement
can be replaced by
the conditions stated in Lemma~\ref{Lem-s4-1}.
\end{remark}

Using the above lemma, we prove the following uniform modulus of
continuity theorem. For convenience, let $B_E(r):=\{x\in\R^N\dvt
\tau
_{E}(x)\leq r\}$
and $B(r)=\{x\in\R^N\dvt \|x\|\leq r\}$ for all $r\geq0$, and $I:=[0,1]^N$.

%
\begin{theorem}\label{Th:UMC}
Let $X=\{X(t), t\in\mathrm{R}^N\}$ be a centered, real-valued Gaussian
random field defined as in (\ref{Hfield}). Then
%
%
\begin{equation}
\label{Eq:M-entropya} \lim_{r \to0} \mathop{\sup_{ s,t \in I}}_{ \tau_{E}(s-t) \le r}
\frac{|X(s) - X(t)|} {\tau_{E}(s-t)\sqrt{\log( 1 +
\tau_{E}(s-t)^{-1})}} = C_{6} \qquad\mbox{a.s.},
\end{equation}
where $C_{6}$ is a positive and finite constant.
\end{theorem}

\begin{pf}
Note that due to monotonicity the limit in the left-hand
side of (\ref{Eq:M-entropya})
exists almost surely, and the key point is that this limit is a
positive and finite constant.

For $t, t'\in I$, let $\beta(t,t')=\tau_{E}(t-t')\sqrt{\log(1+
\tau_{E}(t-t')^{-1})}$
and let
\[
\mathcal{J}(r)=\mathop{\sup_{ t,t'\in I}}_{\tau_{E}(t-t')\le
r}
\frac
{|X(t)-X(t')|} {\beta(t,t')}.
\]

First, we prove that $\lim_{r \to0} \mathcal{J}(r) \le C < \infty$
almost surely.
We introduce an auxiliary Gaussian field:
\[
Y=\bigl\{Y(t,s), t\in I, s\in B_E(r)\bigr\}
\]
defined by $Y(t,s)=X(t+s)-X(t)$, where $r$ is sufficiently
small such that $B_E(r) \subseteq[-1, 1]^N $. Since $X$ has stationary
increments and $X(0)=0$, $d_{X}(s,t)=d_{X}(0,t-s)$
for any $s, t\in\R^N$, the canonical metric $d_Y$ on $T:=I\times
B_E(r)$ associated with $Y$ satisfies the following inequality:
\[
d_Y\bigl((t,s),\bigl(t',s'\bigr)\bigr)
\leq C\min\bigl\{d_{X}(0,s)+d_{X}\bigl(0,s'
\bigr), d_{X}\bigl(s,s'\bigr)+d_{X}
\bigl(t,t'\bigr)\bigr\}
\]
for some constant $C$. Denote the diameter of $T$ in the metric $d_Y$
by $D$.
By Theorem~\ref{Thm:SLND}(i), we have that
\[
D\leq C\mathop{\sup_{s\in B_E(r)}}_{s'\in
B_E(r)}
\bigl(d_{X}(0,s)+d_{X}\bigl(0,s'\bigr)\bigr)
\leq C r
\]
for some constant $C$. Note that by Theorem~\ref{Thm:SLND}(i) and
(\ref{s2-1}), for a
given small $\delta>0$, there is $C>0$ such that
\[
d_{X}(s,t)\leq C\|t-s\|^{1/(a_p+\delta)}.
\]
Therefore, there exists $C>0$ such that for small $\ep>0$, if
$\|t-t'\|< C\ep^{a_p+\delta}$ and $\|s-s'\|<
C\ep^{a_p+\delta}$, then
\[
\bigl(t',s'\bigr)\in\mathrm{O}_{d_Y}
\bigl((t,s),\ep\bigr)=\bigl\{(u,v)\dvt d_{Y}\bigl((u,v),(t,s)\bigr)<
\ep\bigr\}.
\]
Hence, $N_d(T,\ep)$, the smallest number of open $d_Y$-balls of
radius $\ep$ needed to cover $T$, satisfies
\[
N_d(T,\ep)\leq C\ep^{-2N(a_p+\delta)},
\]
for some constant $C>0$. Then one can verify that for some
constant $C>0$
\[
\int_0^D\sqrt{\ln N_d(T, \ep)}
\,\mathrm{d}\ep\leq C r\sqrt{\log \bigl(1+r^{-1}\bigr)}.
\]
It follows from Lemma~2.1 in Talagrand \cite{T95} that for all $u\geq
2Cr\sqrt{\log(1+r^{-1})}$,
\[
\P \Bigl(\sup_{(t,s)\in T}\bigl |X(t+s)-X(t)\bigr |\geq u \Bigr)\leq\exp
\biggl(-\frac{u^2}{4D^2} \biggr).
\]
By a standard Borel--Cantelli argument, we have that for some
positive constant $C<\infty$,
\[
\limsup_{r \to0} \mathop{\sup_{t\in I}}_{\tau_{E}(s-t) \le r}
\frac{|X(s) - X(t)|} {r\sqrt{\log( 1 + r^{-1})}} \le C \qquad\mbox{a.s.}
\]
The monotonicity of the functions $r \mapsto
r \sqrt{\log(1 + r^{-1})}$ implies that $\lim_{r \to0} \mathcal
{J}(r) \le C $ almost surely. Hence, by
Lemma~\ref{Lem-s4-1}, we see that (\ref{Eq:M-entropya}) holds for
a constant $C_6 \in[0, \infty)$.

In order to show $C_6 >0$ it is sufficient to prove that
%
%
\begin{equation}
\label{Eq:UMC-low} \lim_{r\to0+}\mathcal{J}(r)\ge C_7,
\qquad\mbox{a.s.},
\end{equation}
where $C_7 =C_{4}^{-1}\sqrt{2C_5 a_1}$. Recall that $a_1$ is the real
part of
eigenvalue $\lambda_1$. For any $k\geq1$, we let
\[
x_i^{(k)}=\bigl\langle i2^{-n_k}\bigr\rangle,
\qquad i=0,1,2,\ldots, 2^{n_k}
\]
and $r_k=\tau_{E}(\langle2^{-n_k} \rangle)$,
where the sequence $\{n_k\}$ is taken as in Lemma~\ref{lem-s2-1}. Since $0<\tau_{E}(\langle2^{-n_k}\rangle)\to0$ as
$k\to\infty$,
the monotonicity of $\mathcal{J}(r)$ implies that
%
%
\begin{eqnarray}
\lim_{r\to0+} \mathcal{J}(r)&=&\lim_{k\to\infty}\sup
_{s,t\in I,
\tau_{E}(s -t)\leq
r_k}\frac{|X(s)-X(t)|}{\beta(s,t)}
\nonumber
\\
&\geq&\liminf_{k\to\infty}\max_{0\leq i\leq
({1}/{2}) (2^{n_k}-1)}
\frac
{|X(x_{2i+1}^{(k)})-X(x_{2i}^{(k)})|}{r_k\sqrt{\log
(1+r_k^{-1})}}
\\
&=:&\liminf_{k\to\infty}\mathcal{J}_k.
\nonumber
\end{eqnarray}
For any small $\delta>0$, denote
$ C_{8}=C_4^{-1}\sqrt{2C_5(a_1-\delta)}$.
For any $\mu\in(0,1)$, we write
%
%
\begin{eqnarray}
\label{thm-4} &&\P \bigl(\mathcal{J}_k\leq(1-\mu)C_8
\bigr)
\nonumber
\\
&&\quad= \P \biggl( \biggl\{\frac{|X(\langle1-2^{-n_k}\rangle)-X(\langle
1-2^{-n_k+1}\rangle)|}{r_k\sqrt{\log(1+r_k^{-1})}}\leq(1-\mu )C_8
\biggr\}
\\
&&\qquad{} \cap \biggl\{\max_{0\leq
i\leq({1}/{2}) (2^{n_k}-1)-1}\frac{|X(\langle(2i+1)2^{-n_k}\rangle
)-X(\langle(2i)2^{-n_k}\rangle)|}{r_k\sqrt{\log(1+r_k^{-1})}} \leq(1-
\mu)C_8 \biggr\} \biggr).
\nonumber
\end{eqnarray}
Let
%
%
\begin{equation}
\label{thm-3} P_1(k)=\P \biggl(\max_{0\leq
i\leq({1}/{2}) (2^{n_k}-1)-1}
\frac{|X(\langle(2i+1)2^{-n_k}\rangle
)-X(\langle(2i)2^{-n_k}\rangle)|}{r_{k}\sqrt{\log(1+r_{k}^{-1})}} \leq(1-\mu)C_8 \biggr)
\end{equation}
and
%
%
\begin{eqnarray}
\label{thm-2} P_2(k)&=& \P \biggl(\frac{|X(\langle1-2^{-n_k}\rangle)-X(\langle1-2^{-n_k+1}
\rangle)|}{r_{k}\sqrt{\log(1+r_{k}^{-1})}} \leq(1-
\mu)C_8 \Big|X\bigl(\bigl\langle1-2^{-n_k+1}\bigr\rangle\bigr);
\nonumber
\\[-8pt]
\\[-8pt]
&&\phantom{\P \biggl(} X\bigl(\bigl\langle(2i+1)2^{-n_k}\bigr\rangle
\bigr), X\bigl(\bigl\langle (2i)2^{-n_k}\bigr\rangle\bigr), 0\leq i\leq
\frac{1} 2 \bigl(2^{n_k}-1\bigr) -1 \biggr).
\nonumber
\end{eqnarray}
It follows from Theorem~\ref{Thm:SLND} and Lemma~\ref{lem-s2-1} that
\begin{eqnarray*}
&&\operatorname{Var} \biggl(X\bigl(\bigl\langle1-2^{-n_k}\bigr\rangle
\bigr)-X\bigl(\bigl\langle1-2^{-n_k+1} \bigr\rangle\bigr) \big|X\bigl(\bigl
\langle1-2^{-n_k+1}\bigr\rangle\bigr);
\\
&&\qquad{} X\bigl(\bigl\langle(2i+1)2^{-n_k}\bigr\rangle\bigr), X\bigl(
\bigl\langle (2i)2^{-n_k}\bigr\rangle\bigr), 0\leq i\leq\frac{1} 2
\bigl(2^{n_k}-1\bigr) -1 \biggr)
\\
&&\quad\geq C_5\min_{1\leq i\leq
2^{n_k}}\tau_{E}^2
\bigl(\bigl\langle i2^{-n_k}\bigr\rangle\bigr)\geq C_5
C_4^{-2}\tau_{E}^2\bigl(\bigl
\langle2^{-n_k}\bigr\rangle\bigr)=C_5C_4^{-2}
r_k^2.
\end{eqnarray*}
Thus by the fact that the conditional distributions of the
Gaussian process is almost surely Gaussian, and by
Anderson's inequality (see Anderson \cite{A55}) and the definition of $C_8$,
we obtain
\[
P_2(k)\leq\P \Bigl(N(0,1)\leq(1-\mu)\sqrt{ 2
(a_1-\delta)\log \bigl(1+r_k^{-1}\bigr)}
\Bigr),
\]
where $N(0,1)$ denotes a standard normal random variable.
By using the following well-known inequality
\[
(2\pi)^{-{1}/{2}}\bigl(1-x^{-2}\bigr)x^{-1}
\mathrm{e}^{-{x^2}/{2}}\leq\P\bigl(N(0,1)>x\bigr) \leq(2\pi)^{-{1}/{2}}x^{-1}
\mathrm{e}^{-{x^2}/{2}},\qquad \forall x>0,
\]
%
we derive that for all $k$ large enough
%
%
\begin{eqnarray}
\label{thm-6} P_2(k)&\leq& 1-\P \Bigl(N(0,1)> (1-\mu)\sqrt{
2 (a_1 - \delta)\log\bigl(1+r_k^{-1}\bigr)}
\Bigr)
\nonumber
\\[-8pt]
\\[-8pt]
&\leq&1-r_k^{(1-\mu/2)^2(a_1 -\delta)} \leq\exp \bigl(-r_k^{ (1-\mu/2)^2(a_1 -\delta)}
\bigr).
\nonumber
\end{eqnarray}
Combining (\ref{thm-4}) with (\ref{thm-3}), (\ref{thm-2}) and (\ref
{thm-6}), we have that
\[
\P\bigl(\mathcal{J}_k\leq(1-\mu)C_{8}\bigr)
\leq\exp \bigl(-r_k^{(1-\mu
/2)^2(a_1 -\delta)} \bigr) P_1(k).
\]
%
By repeating the above argument, we obtain
%
%
\begin{equation}
\label{thm-7} \P \bigl(\mathcal{J}_k\leq(1-\mu)C_{8}
\bigr)\leq \exp \biggl(-\frac{2^{n_k}-1}{2} r_k^{ (1-\mu/2)^2(a_1 -\delta
)} \biggr)
\leq\exp \bigl(-C 2^{\mu n_k/2} \bigr),
\end{equation}
where the last inequality follows from the estimate:
\[
r_k^2 =\tau_{E}^2\bigl(\bigl
\langle2^{-n_k}\bigr\rangle\bigr)\geq C_{1}^2\bigl \|
\bigl\langle2^{-n_k}\bigr\rangle\bigr \|^{{2}/{(a_1-\delta)}}\geq C_{1}^22^{{-2n_k}/{(a_1-\delta)}}.
\]
By (\ref{thm-7}) and the Borel--Cantelli lemma, we have
$\liminf_{k\to\infty}\mathcal{J}_k\geq(1-\mu)C_{8}$ {a.s.}
Letting $\mu\to0$ and $\delta\to0$ yields (\ref{Eq:UMC-low}).
The proof of Theorem~\ref{Th:UMC} is completed.
\end{pf}

\section{Laws of the iterated logarithm}\label{sec5}

For any fixed $t_0 \in\R^N$ and a family of neighborhoods $\{
{O}(r)\dvt
r>0\}$ of $ 0 \in\R^N$ whose diameters go to 0
as $ r \to0$, we consider in this section the corresponding
local modulus of continuity of $X$ at $t_0$
\[
\omega(t_0, r) = \sup_{s \in {O}(r)}\bigl |X(t_0
+s)-X(t_0)\bigr |.
\]
Since $X$ is anisotropic, the rate at
which $\omega(t_0, r)$ goes to 0 as $r \to0$ depends on
the shape of ${O}(r)$. A natural choice of ${O}(r)$ is $B_E(r)$.

For specification and simplification, in this section, let $E$ be a
Jordan canonical form of (\ref{Eq:J0}), which satisfies all
assumptions in Section~\ref{sec2}. Recall that $\tilde{l}_j$ is the
size of
$J_j$. For any $i=1,2,\ldots,N$, if
$\tilde{l}_1+\cdots+\tilde{l}_{j-1}+1 \leq i\leq
\tilde{l}_1+\cdots+\tilde{l}_{j}$, then
\[
e_i=\{\underbrace{0,\ldots,0,1}_i,0, \ldots,0\}\in
W_j.
\]
%

The following theorem characterizes the exact local modulus of continuity
of $X$.

%
\begin{theorem}\label{Th:LIL2}
There is a positive and finite constant $C_{9}$ such that for
every $t_0\in\R^N$ we have
%
%
\begin{equation}
\label{Eq:LIL2} \lim_{r\to0+} \sup_{{s-t_0\in B_E(r)}}
\frac{|X(s)-X(t_0)|}{\tau_E(s-t_0)\sqrt{\log\log(1 +
\tau_E(s-t_0)^{-1})}} =C_{9} \qquad\mbox{a.s.}
\end{equation}
%
\end{theorem}

In order to show this result, we will make use of the following
lemmas.

%
\begin{lemma}\label{Lem:Tail}
There exist positive and finite constants $u_0$ and
$C_{10}$ such that for all $t_0 \in\R^N$, $u \ge u_0$ and
sufficiently small $r>0$,
\[
\P \Bigl(\sup_{s\in B_E(r)} \bigl |X(t_0+s)-X(t_0)\bigr |
\ge u r\sqrt{\log\log \bigl(1+ r^{-1}\bigr)} \Bigr) \le
\mathrm{e}^{-C_{10} u^2\log\log(1+r^{-1})}.
\]
\end{lemma}

\begin{pf}
We introduce an auxiliary Gaussian field
$Y=\{Y(s), s\in B_E(r)\}$
defined by $Y(s)=X(t_0+s)-X(t_0)$. Since $X$ has
stationary increments and $X(0)=0$,
we have $d_Y(s,s')=d_{X}(s,s')$ for all $s, s'\in\R^N$.
Denote the diameter of $B_E(r)$ in the metric $d_Y$ by $D$. It follows from
Lemma~\ref{Lem:Var1} that $D \leq Cr$
for some finite constant $C$. Note that the decomposition of
$x=(x_1,x_2,\ldots,x_N)
\in B_{E}(r)$ in $W_j$ is
\[
\bar{x}_j=(0,\ldots,0,x_{\tilde{l}_1+\cdots+\tilde
{l}_{j-1}+1},\ldots,x_{\tilde{l}_1+\cdots+\tilde{l}_{j}},
0,\ldots,0).
\]
For any $j=1,2,\ldots,p$, let $l_j=\tilde{l}_j$ if $J_j$ is a Jordan
cell matrix as in (\ref{Eq:J1}) or $l_j=\tilde{l}_j/2$
if $J_j$ is of the form (\ref{Eq:J2}). By Lemma~\ref{lem-s2-3} and
(\ref{s2-3}), we have that for sufficiently
small $r$,
\[
\|\bar{x}_j \|^{1/a_j} \bigl |\ln\|\bar{x}_j \|
\bigr |^{-(l_j-1)/a_j}\leq C r.
\]
This implies that there exists a constant $C$, which may depend on
$a_j$, such
that for all $i$ with $\tilde{l}_1+\cdots+\tilde{l}_{j-1}+1\leq
i\leq\tilde{l}_1+\cdots+\tilde{l}_{j}$,
\[
|x_i|\leq C r^{a_j}|\ln r|^{ l_j-1 }.
\]
Therefore $B_E(r)\subset[-h, h]$ for sufficiently small $r>0$,
where $h=(h_1,h_2,\ldots, h_N)$ with $h_i=Cr^{a_j}|\ln
r|^{ l_j-1}$ as $\tilde{l}_1+\cdots+\tilde{l}_{j-1}+1\leq i\leq
\tilde{l}_1+\cdots+\tilde{l}_{j}$. Furthermore,
from (\ref{s2-3}), we have that for any $x=(x_1,x_2,\ldots,x_N)\in\R
^N$ and
sufficiently small $\ep>0$, if
\[
|x_i|< \biggl(\frac{\ep}{N\mu} \biggr)^{a_j}\biggl |\ln
\biggl(\frac{\ep
}{N\mu}\biggr)\biggr |^{-(l_j-1)/a_j},
\]
for $\tilde{l}_1+\cdots+\tilde{l}_{j-1}+1\leq i\leq\tilde
{l}_1+\cdots+\tilde{l}_{j}$, where $\mu$ is a constant
whose value will be determined later, then
%
%
\begin{eqnarray}
\label{Eq:tau4} %
\tau_{E}(\vec{x}_i) &\le& C
\frac{\ep}{N\mu} \biggl|\ln\biggl(\frac{\ep
}{N\mu}\biggr)\biggr |^{-(l_j-1)} \biggl|\ln
\biggl[ \biggl(\frac{\ep}{N\mu} \biggr)^{a_j}\biggl |\ln\biggl(
\frac
{\ep}{N\mu}\biggr) \biggr|^{-(l_j-1)/a_j} \biggr] \biggr|^{(l_j-1)/a_j}
\nonumber
\\[-8pt]
\\[-8pt]
&\leq& C \frac{\ep}{N\mu},
\nonumber
\end{eqnarray}
where\vspace*{-5pt} $\vec{x}_i=(\underbrace{0,\ldots,0,
x_i}_i,0,\ldots,0)\in\R^N$. Then by (\ref{s2-6}) and (\ref
{Eq:tau4}), there exists a
constant $C>0$ such that
\[
\tau_{E}(x)=\tau_{E} \Biggl(\sum
_{i=1}^N\vec{x}_i \Biggr)\le C \sum
_{i=1}^N \frac{\ep}{N\mu} \leq C
\frac{\ep} \mu.
\]
By using Lemma~\ref{Lem:Var1} again, we have
\[
d_Y(0,x) \le C \tau_{E}(x) \leq C_{11}\ep/
\mu.
\]
Now we take $\mu> C_{11}$, then $x\in\mathrm{O}_{d_Y}(\ep)$ implies
$[0,x)\subset\mathrm{O}_{d_Y}(\ep)$.
Therefore the smallest number of open $d_Y$-balls of
radius $\ep$ needed to cover $B_E(r):=T$, denoted by $N_d(T,\ep)$, satisfies
\[
N_d(T,\ep)\leq C\prod_{j=1}^p
\biggl(\frac{r^{a_j}}{({\ep}/{(\mu N)})^{a_j}} \biggl|\ln\biggl(\frac{\ep}{\mu N}\biggr)\biggr |^{(l_j-1)/a_j} |
\ln r |^{(l_j-1)} \biggr)^{l_j}
\]
for some constant $C>0$. Then one can verify that
\[
\int_0^D\sqrt{\ln N_d(T, \ep)}
\,\mathrm{d}\ep\leq C r\sqrt{\log \log\bigl(1+r^{-1}\bigr)}.
\]
It follows from \cite{T95}, Lemma~2.1, that for all sufficiently large $u$,\vspace*{1pt}
\begin{eqnarray*}
&&\P \Bigl(\sup_{s\in B_E(r)}\bigl |X(t_0+s)-X(t_0)\bigr |
\geq ur\sqrt{\log\log \bigl(1+r^{-1}\bigr)} \Bigr)
\\
&&\quad\leq\exp \bigl(- C_{10}u^2\log\log
\bigl(1+r^{-1}\bigr) \bigr).
\end{eqnarray*}
This finishes the proof of Lemma~\ref{Lem:Tail}.
\end{pf}

%
\begin{lemma}\label{lem-5-2}
There is a constant $C_{12} \in[0, \infty)$ such that for every fixed
$t_0 \in\R^N$,\vspace*{1pt}
%
%
\begin{equation}
\label{Eq:LIL07a} \lim_{\ep\to0} \sup_{s-t_0\in B_E(\ep) }
\frac{|X(s) - X(t_0)|} {\tau_E(s- t_0)\sqrt{\log\log( 1 +
\tau_E(s-t_0)^{-1})}} = C_{12} \qquad\mbox{a.s.}
\end{equation}
\end{lemma}

\begin{pf}
By Lemma~\ref{Lem-s4-1},
it is sufficient to prove
%
%
\begin{equation}
\label{s5-1} \lim_{\ep\to0}\sup_{s-t_0\in B_E(\ep)}
\frac{|X(s) - X(t_0)|} {\tau_E(s- t_0)\sqrt{\log\log(
1 + \tau_E(s- t_0)^{-1})}}\leq C<\infty,
\end{equation}
for some constant $C>0$. Let $\ep_n=\mathrm{e}^{-n}$, consider the event
\[
E_n= \biggl\{\sup_{ s-t_0\in B_E(\ep)} \frac{|X(s) - X(t_0)|} {\ep_n\sqrt{\log\log(
1 + \ep_n^{-1})}}> u
\biggr\},
\]
where $u>C_{10}^{-1/2}$ is a constant. By Lemma~\ref{Lem:Tail}, we have
$\P(E_n)\leq\mathrm{e}^{-C_{10}u^2\log n}$
for all sufficiently large $n$. Hence, the Borel--Cantelli lemma implies
\[
\limsup_{\ep\to0}\sup_{s\in I, s-t_0\in B_E(\ep)}
\frac{|X(s) - X(t_0)|} {\ep\sqrt{\log\log( 1 + \ep^{-1})}}\leq u.
\]
%
This and a monotonicity argument yield (\ref{s5-1}).
\end{pf}

We will also need the following truncation inequalities which
extend a result in Luan and Xiao \cite{LuanXiao10}.

%
\begin{lemma}\label{Lem:Truncation}
For a given $N\times N$ matrix $E$, there
exists a constant $r_0>0$ such that for any $u > 0$
and any $t \in\R^N$ with $\tau_{E}(t) u \le r_0$, we have
%
%
\begin{equation}
\label{Eq:Trun1} \int_{\{\tau_{{E'}}(\xi) < u\}}\langle t,\xi\rangle^2
\frac{\mathrm{d} \xi}{\psi(\xi)^{2+Q}} \le 3 \int_{\R^N} \bigl(1 - \cos\langle t,
\xi\rangle\bigr) \frac{d
\xi}{\psi(\xi)^{2+Q}}.
\end{equation}
\end{lemma}

\begin{pf}
Let $M=\max\{\|x\|, x\in S_{E}\}$,
$K(r)=\max\{\|x\|,\tau_{{E'}}(x)\leq r\}$. Since $S_{E}$ is compact
set without $0$ and $\tau_{{E'}}(\cdot)$ is continuous, $M>0$ and
$K(r)$ continuous with $K(0)=0$, $K(r)\to\infty$ as $r\to\infty$.
Therefore, there exists $r_0>0$ such that $MK(r)\leq1$ for all
$r<r_0$. By using the inequality $u^2\leq3(1-\cos u)$ for all real
numbers $|u|\leq1$, we derive that if $\tau_{E}(t)u\le r_0$, then
\begin{eqnarray*}
&&\int_{\{\tau_{{E'}}(\xi) < u\}}\langle t,\xi\rangle^2
\frac{\mathrm{d} \xi}{\psi(\xi)^{2+Q}}
\\
&&\quad = \int_{\{
\tau
_{{E'}}(\xi) < u\}}\bigl\langle
\tau_{E}^{E}(t)l_{E}(t),\xi\bigr
\rangle^2 \frac{\mathrm{d} \xi}{\psi
(\xi)^{2+Q}}
\\
&&\quad = \int_{\{\tau_{{E'}}(\xi) < u\}}\bigl\langle l_{E}(t),
\tau_{E}^{E'}(t)\xi\bigr\rangle^2
\frac{\mathrm{d} \xi}{\psi
(\xi)^{2+Q}} =\tau_{E}^2(t)\int_{\{\tau_{{E'}}(\xi)<\tau_{E}(t) u\}}
\bigl\langle l_{E}(t),\xi\bigr\rangle^2 \frac{\mathrm{d} \xi}{\psi(\xi)^{2+Q}}
\\
&&\quad \leq 3\tau_{E}^2(t)\int_{\{\tau_{{E'}}(\xi)<\tau_{E}(t) u\}
}
\bigl(1-\cos\bigl\langle l_{E}(t),\xi\bigr\rangle \bigr)
\frac{\mathrm{d} \xi}{\psi(\xi)^{2+Q}}
\\
&&\quad =  3 \int_{\{\tau_{{E'}}(\xi)<u\}} \bigl(1-\cos\bigl\langle
l_{E}(t),\tau_{E}^{E'}(t)\xi\bigr\rangle \bigr)
\frac{\mathrm{d}
\xi}{\psi(\xi)^{2+Q}},
\end{eqnarray*}
which equals
%
%
\begin{equation}
\label{trun2} 3 \int_{\{\tau_{{E'}}(\xi)<u\}} \bigl(1-\cos\langle t,\xi\rangle
\bigr) \frac{\mathrm{d} \xi}{\psi(\xi)^{2+Q}}\leq3 \int_{\R
^N} \bigl(1-\cos\langle
t,\xi\rangle \bigr) \frac{\mathrm{d} \xi}{\psi(\xi)^{2+Q}}.
\end{equation}
The proof of this lemma is complete.
\end{pf}

Now we are ready to prove Theorem~\ref{Th:LIL2}.

\begin{pf*}{Proof of Theorem~\ref{Th:LIL2}}
By Lemma~\ref{lem-5-2} and
the stationary increments property of $X$, it only remains to show\vspace*{1pt}
%
%
\begin{equation}
\label{s5-0-1} \lim_{\ep\to0+} \sup_{ s\in B_E(\ep)}
\frac{|X(s)|}{\tau
_E(s)\sqrt{\log\log(1 +
\tau_E(s)^{-1})}}\geq C
\end{equation}
for some constant $C > 0$.

For any $0<\mu<1$ and $n\geq1$, we define
$s_n=(0,\ldots,0,\mathrm{e}^{-a_pn^{1+\mu}})\in\R^N$.
By (\ref{s2-3})
%
%
\begin{equation}
\label{s5-0} C_{3}^{-1}\mathrm{e}^{-n^{1+\mu}}\bigl |a_pn^{1+\mu
}\bigr |^{-{(l_p-1)}/{a_p}}
\leq\tau _{E}(s_n) 
\leq C_{3}
\mathrm{e}^{-n^{1+\mu}}\bigl |a_pn^{1+\mu}\bigr |^{{(l_p-1)}/{a_p}}.
\end{equation}

For every integer $n\geq1$, let $d_n=\exp(n^{1+\mu}+n^\mu)$. Denote
$U=\exp(\mu(n-1)^\mu)$. Notice that as $n\to\infty$,
\begin{eqnarray*}
\tau_{E}(Us_n)d_{n-1}&\leq& C U^{1/a_p}
\bigl |a_pn^{1+\mu}-\mu(n-1)^\mu\bigr  |^{{(l_p-1)}/{a_p}}\exp
\bigl(-n^{1+\mu}+(n-1)^{1+\mu}+(n-1)^\mu\bigr)
\\
&\leq&C \bigl |a_pn^{1+\mu}-\mu(n-1)^\mu
\bigr |^{{(l_p-1)}/{a_p}}\exp \biggl(-\mu\biggl(1-\frac{1}{a_p}\biggr)
(n-1)^\mu \biggr)\to 0.
\end{eqnarray*}
It follows from Lemma~\ref{Lem:Truncation}, Lemma~\ref{Lem:Var1} and
(\ref{s2-3}) that\vspace*{1pt}
%
%
\begin{eqnarray}
\label{s5-8} 
&& \int_{\{\tau_{{E'}}(\xi)\leq d_{n-1}\}}\langle s_n, \xi
\rangle^2\frac{\mathrm{d}\xi}{\psi(\xi)^{2+Q}}\nonumber
\\
&&\quad = U^{-2}\int
_{\{\tau _{{E'}}(\xi)\leq
d_{n-1}\}}\langle s_n U, \xi\rangle^2
\frac{\mathrm{d}\xi}{\psi
(\xi )^{2+Q}}
\nonumber
\\
&&\quad \leq CU^{-2}d_{X}^2(Us_n,0)\leq
CU^{-2}\tau_{E}^2(Us_n)
\\
&&\quad \leq CU^{-2}U^{{2}/{a_p}} \bigl |\ln\|s_n\|
\bigr |^{2{(l_p-1)}/{a_p}} \bigl |\ln\|Us_n\| \bigr |^{2{(l_p-1)}/{a_p}}\tau
_{E}^2(s_n)\nonumber
\\
&&\quad \leq C \exp \biggl(-\biggl(1-\frac{1}{a_p}\biggr)\mu(n-1)^\mu
\biggr)\tau _{E}^2(s_n)
\nonumber
\end{eqnarray}
for $n$ large enough.
On the other hand, noting that $\psi$ is
$E'$-homogeneous, by using \cite{BMS07}, Proposition~2.3, we
obtain that
\[
\int_{\{\tau_{{E'}}(\xi)>d_n\}}\frac{\mathrm{d}\xi}{\psi(\xi
)^{2+Q}}=\int_{d_n}^\infty
\,\mathrm{d}r \int_{S_{E'}}\frac{1}{r^{2+Q}\psi(\theta)}r^{Q-1}
\sigma(\mathrm{d} \theta)\leq Cd_n^{-2},
\]
since $\sigma(\mathrm{d}\theta)$ is a finite measure on $S_{E'}$.
Furthermore,
\begin{eqnarray*}
d_n^{-2}&=&\mathrm{e}^{-2n^{1+\mu}-2n^{\mu}}=\mathrm{e}^{-2n^{1+\mu
}}\bigl |
\ln\mathrm{e}^{-a_pn^{1+\mu}}\bigr |^{-2{(l_p-1)}/{a_p}}\bigl |a_pn^{1+\mu
}\bigr |^{2{(l_p-1)}/{a_p}}
\mathrm{e}^{-2n^{\mu
}}
\\
&\leq& C\tau_{E}^2(s_n)\bigl |n^{1+\mu}\bigr |^{2{(l_p-1)}/{a_p}}
\mathrm{e}^{-2n^\mu},
\end{eqnarray*}
when $n$ is large enough. Therefore, for sufficiently large $n$,
%
%
\begin{equation}
\label{s5-9} \int_{\{\tau_{{E'}}(\xi)>d_n\}}\frac{\mathrm{d}\xi}{\psi(\xi
)^{2+Q}}\leq C
\tau_{E}^2(s_n)\mathrm{e}^{-n^\mu}.
\end{equation}
Now we decompose $X$ into two independent parts as follows.
%
%
\begin{equation}
\label{s5-4} \widetilde{X}_n(t)=\int_{\{\tau_{{E'}}(\xi)\notin(d_{n-1},d_n]\}
}
\bigl(\mathrm{e}^{\mathrm{i}\langle t,\xi\rangle}-1 \bigr) \frac{\widetilde{\mathcal{M}} (\mathrm{d}\xi)} { \psi(\xi)^{1 + Q/2}}
\end{equation}
and
%
%
\begin{equation}
\label{s5-5} X_n(t)=\int_{\{\tau_{{E'}}(\xi)\in(d_{n-1},d_n]\}} \bigl(
\mathrm{e}^{\mathrm{i}\langle t,\xi\rangle}-1 \bigr) \frac{\widetilde{\mathcal{M}} (\mathrm{d}\xi)} { \psi(\xi)^{1 + Q/2}}.
\end{equation}
Notice that the random fields $\{X_n(t), t \in\R^N\}$, $n=1,2, \ldots
$ are independent.

Let
\[
I_1(n)=\frac{|X_n(s_n)|}{\tau_{E}(s_n)\sqrt{\log\log(1+\tau
_{E}(s_n)^{-1})}}
\]
and
\[
I_2(n)=\frac{|\widetilde{X}_n(s_n)|}{\tau_{E}(s_n)\sqrt{\log\log
(1+ \tau_{E}(s_n)^{-1})}}.
\]
Then
%
%
\begin{eqnarray}
\label{s5-6} \lim_{\ep\to0+} \sup_{ s\in B_E(\ep)}
\frac{|X(s)|}{\tau_E(s)\sqrt{\log\log(1 +
\tau_E(s)^{-1})}} &\geq&\limsup_{n\to\infty} \frac{|X(s_n)|}{\tau_{E}(s_n)\sqrt{\log\log(1+\tau
_{E}(s_n)^{-1})}}\qquad
\nonumber
\\[-8pt]
\\[-8pt]
&\geq&\limsup_{n\to\infty}I_1(n)-\limsup
_{n\to
\infty}I_2(n).
\nonumber
\end{eqnarray}
By using (\ref{s5-8}), (\ref{s5-9}) and the same argument in the
proof of Theorem~5.5 in
\cite{MWX10}, we can readily get that
%
\begin{equation}
\label{s5-10} \limsup_{n\to\infty}I_2(n)=0, \qquad
\mbox{a.s.}
\end{equation}
In order to estimate $\limsup_{n\to\infty}I_1(n)$, using
Lemma~\ref{Lem:Var1} again, we have that
\[
\E \bigl(X_n(s_n) \bigr)^2\leq
d_{X}^2(s_n,0)\leq C_{13}
\tau_{E}^2(s_n).
\]
Again, by the corresponding argument in the proof of Theorem~5.5 in
\cite{MWX10}, it is
easy to get that
%
%
\begin{equation}
\label{s5-11} \limsup_{n\to\infty}I_1(n)\geq\sqrt{2
C_{13}} \qquad\mbox{a.s.}
\end{equation}
Hence, (\ref{s5-0-1}) follows from
(\ref{s5-6}), (\ref{s5-10}) and (\ref{s5-11}).
\end{pf*}

\section{Examples}\label{sec6}

Finally, we provide two examples of operator scaling Gaussian
random fields with stationary increments to illustrate our results
and compare them with those in Meerschaert, Wang and Xiao \cite{MWX10}.
In particular, Example~\ref{exa6.2}
shows that the regularity properties of $X$ depend
subtly on its scaling exponent $E$.

\begin{Example}\label{exa6.1}
If $E$ has a Jordan canonical form (\ref{Eq:J0}) such that, for all
$k=1,2,\ldots, p$, $\tilde{l}_k=1$ if $J_k$ is a Jordan cell matrix
and $\tilde{l}_k=2$ if $J_k$ is not a Jordan cell matrix. Then for
any $t=(t_1,\ldots, t_N)\in\R^N$, by (\ref{s2-3}) and Lemma~\ref{lem-s2b}, we have
\[
\tau_{E}(t) \asymp\sum_{i=1}^N
|t_i|^{1/a_i},
\]
where $a_i$ is the real part of eigenvalue(s) corresponding to $J_k$
such that $\sum_{j=1}^{k-1}\tilde{l}_j+1\leq i\leq
\sum_{j=1}^{k}\tilde{l}_j$. Therefore, in this case, Theorems~\ref
{Th:UMC} and~\ref{Th:LIL2} are of the same form as the corresponding results in
Meerschaert, Wang and Xiao \cite{MWX10}.
\end{Example}

\begin{Example}\label{exa6.2} We consider the Gaussian random field
$\{X(t), t \in\R^2\}$ defined by (\ref{Hfield}) with scaling
exponent $E$, a Jordan
matrix, as follows
\[
E=\lleft( %
\begin{array} {@{}c@{\quad}c@{}} a &0
\\
1 &a \end{array} %
\rright),
\]
where $a>1$ is a constant. Then $p=1$ and $\tilde{l}_1 = 2$. For any
$t > 0$,
by straightforward computations, we have
\[
t^E=t^a\lleft( %
\begin{array} {@{}c@{
\quad}c@{}} 1 &0
\\
\ln t &1 \end{array} %
\rright).
\]
According to Lemma~6.1.5 in \cite{RVbook}, the norm $\|\cdot\|_E$
induced by $E$ is defined as that for any $x\in\R^2$
\[
\|x\|_E=\int_0^1
\frac{\|t^Ex\|}{t}\,\mathrm{d} t.
\]
Note that we can uniquely represent $x\in\R^2$ as
$(0,s)$ or $(s,\theta s)$ for some $s\in\R, \theta\in\R$. %
When $x=(0,s)$,
%
%
\begin{equation}
\label{ex2-1} \|x\|_E=\int_0^1
|s|t^{a-1}\,\mathrm{d}t=\frac{|s|}{a},
\end{equation}
and when $x=(s,\theta s)$,
%
%
\begin{equation}
\label{ex2-2} \|x\|_E=\int_0^1
|s|t^{a-1}\sqrt{1+(\theta+\ln t)^2}\,\mathrm{d} t=:|s|
\alpha(\theta).
\end{equation}
It is easy to see that $\alpha(\theta)$ is continuous on
$\theta\in\R$ with $\alpha(\theta)>1/a$ and that
$|\theta|/\alpha(\theta)$ is bounded since $|\theta|/\alpha(\theta)$
is continuous and
%
%
\begin{equation}
\label{ex2-6} \lim_{\theta\to\infty}\frac{|\theta|}{\alpha(\theta)}=a.
\end{equation}
We have
$\alpha:=\min_{\theta}\alpha(\theta)>1/a$. From (\ref{ex2-1}) and
(\ref{ex2-2}), we have
\[
S_E=\bigl \{x\dvt\|x\|_E=1\bigr \}= \biggl\{\pm\binom{0} {a},
\pm\frac{1}{\alpha(\theta)} \binom{1} {\theta} 
\dvt \theta\in\R \biggr\},
\]
and $\R^2=\{s^E y\dvt s\geq0, y\in S_E\}$.

To unify the notation, we set
\[
\frac{\theta} {\alpha(\theta)}=\pm a \quad\mbox{and}\quad\frac{1}{\alpha
(\theta)}=0
\]
when $\theta=\pm\infty$. Then for any $x\in\R^2$
with $\tau_{E}(x)=s$, there exists $\theta\in[-\infty,+\infty]$ such
that
%
%
\begin{equation}
\label{ex2-3} x=\pm s^E \frac{1}{\alpha(\theta)} \binom1 \theta
=\pm\frac{s^a}{\alpha(\theta)} \binom1 {\theta+\ln s},
\end{equation}
where $s^a\ln s|_{s=0}:=0$ and the sign $+$ or $-$ depends on $x$.

Now we reformulate Theorem~\ref{Th:UMC} and Theorem~\ref{Th:LIL2} for
the present case. For convenience, we express the vector $y \in\R^2$
in terms
of $s=\tau_E(y)$ and $\theta$ by
\[
y=y(s,\theta, w)=(-1)^w \biggl(\frac{s^a}{\alpha(\theta)},
\frac
{s^a}{\alpha(\theta)}(\theta+\ln s) \biggr),
\]
where $w\in\{0, 1\}$.
\end{Example}

\renewcommand{\theConclusion}{A}
\begin{Conclusion}\label{conA} Let $I=[0, 1]^2$. Then
%
%
\begin{equation}
\label{ex2-4} \lim_{r\to0+}\mathop{\sup_{s\le r, \theta\in[-\infty,+\infty
]}}_{w\in\{0,
1\}, x,x+y\in I}
\frac{|X(x+y(s,\theta, w))-X(x)|} {
s\sqrt{\log(1+s^{-1})}}=C_{17}\qquad\mbox{a.s.},
\end{equation}
and that for any $x_0\in I$,
%
%
\begin{equation}
\label{ex2-5} \lim_{r\to0+} \sup_{s\leq r, \theta\in[-\infty,+\infty],
w\in\{0,1\}}
\frac{|X(x_0+y(s,\theta, w))-X(x_0)|}{s\sqrt{\log\log
(1 + s^{-1})}} =C_{18} \qquad\mbox{a.s.},
\end{equation}
where $C_{17}$ and $C_{18}$ are positive and finite constants.
\end{Conclusion}

Next we describe the asymptotic behavior of $\tau_E(y)$ as $\|y\| \to
0$ along three types of curves in $\R^2$:
\begin{enumerate}[(iii)]
\item[(i)] If $\theta=-\ln s+c$ for a constant $c\in\R$, then
$y=y(s, \theta,
w)=(-1)^w(s^a/\alpha(\theta), cs^a/\alpha(\theta))$ satisfies
\[
\|y\|=\frac{\sqrt{1+c^2}s^a}{\alpha(\theta)}=\frac{\sqrt {1+c^2}s^a}{\alpha(c-\ln s)}.
\]
This, together with (\ref{ex2-6}), implies that as $\|y\|\to0$,
%
%
\begin{equation}
\label{ex2-9} s=\tau_E(y)\sim\|y\|^{1/a} \bigl |\ln\|y\|
\bigr |^{1/a},
\end{equation}
where the notation ``$\sim$'' means that as $\|y\|\to0$ the quotient
of the two sides of $\sim$ goes to a positive constant.
\item[(ii)] If $\theta=\pm\infty$, then $y(s, \theta, w)=(-1)^w
(0, a s^a)$ and
%
%
\begin{equation}
\label{ex2-8} s=\tau_E(y)=\frac{1}{a^{1/a}}\|y\|^{1/a}.
\end{equation}
\item[(iii)] If $\theta$ is fixed in $(-\infty, +\infty)$,
then for $y=y(s,\theta, w)$,
\[
\|y\|=\frac{s^a}{\alpha(\theta)}\sqrt{1+(\theta+\ln s)^2},
\]
which implies that as $\|y\|\to0$,
%
%
\begin{equation}
\label{ex2-7} s=\tau_E(y)\sim\|y\|^{1/a} \bigl |\ln\|y\|
\bigr |^{-1/a}.
\end{equation}
\end{enumerate}

In the following, we derive the exact uniform moduli of continuity
of $X(x)$ by using the norm $\|\cdot\|$ in three different cases
which are intuitively corresponding to the three types mentioned
above. These results illustrate the subtle changes of the regularity
properties of $X$. For the exact local moduli of continuity, similar
results are true as well. In order not to make the paper too
lengthy, we leave it to interested readers.

\renewcommand{\theConclusion}{B}
\begin{Conclusion}\label{conB} (1) If $I_1=\{(t, t)\dvt t\in[0, 1]\}$,
then
%
%
\begin{equation}
\label{ex2-11} \lim_{\|y\|\to0} \sup_{x,x+y\in I_1}
\frac{|X(x+y)-X(x)|} {(\|y\|
|\ln
\|y\| |)^{1/a}\sqrt{\log(1+\|y\|^{-1})}}=C_{19}\in(0, \infty)\qquad\mbox{a.s.}
\end{equation}

(2) If $I_2=\{(0, t)\dvt t\in[0, 1]\}$, then
%
%
\begin{equation}
\label{ex2-12} \lim_{\|y\|\to0} \sup_{x,x+y\in I_2}
\frac{|X(x+y)-X(x)|} {
\|y\|^{1/a}\sqrt{\log(1+\|y\|^{-1})}}=C_{20}\in(0, \infty) \qquad\mbox{a.s.}
\end{equation}

(3) Let $\theta_0\in\operatorname{arg}\min_\theta\alpha(\theta
)=\{
\vartheta, \alpha(\vartheta)\leq\alpha(\theta),\theta\in
[-\infty, +\infty]\}$. Then
%
%
\begin{equation}
\label{ex2-15} \lim_{\|y\|\to0}\mathop{\sup_{y=y(r,\theta_0,0)}}_{
x,x+y\in I}
\frac{ |\ln\|y\| |^{1/a}|X(x+y)-X(x)|} {
\|y\|^{1/a}\sqrt{\log(1+\|y\|^{-1})}}=C_{21}\in(0, \infty)\qquad\mbox{a.s.}
\end{equation}
\end{Conclusion}

\begin{pf} (1) Observe that, in the proof of (\ref{Eq:UMC-low}) in
the case of $N=2$,
one can choose the sequences of $\{x_i^{(k)}\}$ such that all the
points $x_i^{(k)}$ and the differences
$x_{i+1}^{(k)}-x_{i}^{(k)} $ lie in $I_1=\{(t, t)\dvt t\in[0, 1]\}$.
Therefore, the proof of (\ref{Eq:UMC-low})
essentially shows that
%
%
\begin{equation}
\label{ex2-10} \lim_{y\to0}\sup_{x,x+y\in I_1}
\frac{|X(x+y)-X(x)|} {
\tau_E(y)\sqrt{\log(1+\tau_E(y)^{-1})}}\ge C >0 \qquad\mbox{a.s.}
\end{equation}
Thanks to the formula (\ref{ex2-9}), we see that (\ref{ex2-11}) follows
from the proof of Theorem~\ref{Th:UMC}.

(2) To prove (\ref{ex2-12}), choose $x_i^{(n)}=(0, i2^{-n})$ for
$i=0,1,\ldots,2^n$. Then by some obvious modifications, one can
easily check that (\ref{ex2-10}) is also true with $I_2$ instead of
$I_1$. Therefore, by (\ref{ex2-8}), (\ref{ex2-12}) also follows from
in the proof of Theorem~\ref{Th:UMC}.

(3) Note that $\alpha(\theta)$ is continuous and as $\theta\to
\infty$, $\alpha(\theta)/|\theta|\to a$.
The set $\operatorname{arg}\min_\theta\alpha(\theta)$ is not
empty. Let
$\alpha_0=\alpha(\theta_0)$ and
\[
x_i^{(n)}=iy\bigl(2^{-n},\theta_0,0
\bigr)+(0,1)= \biggl(\frac{i2^{-an}}{\alpha_0}, \frac{i2^{-an}}{\alpha_0}( \theta_0-n
\ln 2)+1 \biggr)
\]
for
$i=0,1,2,\ldots,K_n$, where
\[
K_n=\max\bigl\{i, x_i^{(n)}\in[0,
1]^2\bigr\}.
\]
Manifestly, for sufficiently large $n$, $K_n>2^n$. Let
$r_n:=\tau_E(y(2^{-n},\theta_0,0))$. Then
\begin{eqnarray*}
&&\lim_{r\to0}\mathop{\sup_{y=y(r,\theta_0,0)}}_{
x,x+y\in I}
\frac{|X(x+y)-X(x)|} {r\sqrt{\log(1+r^{-1})}}
\\
&& \quad\geq\liminf_{n\to\infty}\max_{0\leq i\leq
K_n-1}
\frac{|X(x_{i+1}^{(n)})-X(x_{i}^{(n)})|}{r_n\sqrt{\log(1+r_n^{-1})}} =:\liminf_{k\to\infty}\mathcal{J}_n.
\end{eqnarray*}
Note that for $k\geq1$,
\[
ky\bigl(2^{-n},\theta_0,0\bigr)= \biggl(
\frac{k2^{-an}}{\alpha_0}, \frac{k2^{-an}}{\alpha_0}( \theta_0-n\ln 2) \biggr).
\]
There exist some $\theta\in(-\infty, \infty)$, $w\in\{0, 1\}$ and
$s=\tau_E(ky(2^{-n},\theta_0,0))$ such that
\[
ky\bigl(2^{-n},\theta_0,0\bigr)=y(s,\theta, w),
\]
which implies that $w=0$ and
\[
\frac{s^a}{\alpha(\theta)}=\frac{k2^{-an}}{\alpha_0}.
\]
Because $\alpha_0=\min_\theta\alpha(\theta)$
\[
s=\tau_E\bigl(ky\bigl(2^{-n},\theta_0,0\bigr)
\bigr)\geq 2^{-n} :=r_n.
\]
Therefore, from
Theorem~\ref{Thm:SLND} and Lemma~\ref{lem-s2-1}, we obtain that
\[
\operatorname{Var} \bigl(X\bigl(x_{i+1}^{(n)}
\bigr)-X\bigl(x_i^{(n)}\bigr) |X\bigl(x_k^{(n)}
\bigr), 0\leq k\leq i \bigr) \geq C_5\min_{1\leq k\leq
i+1}
\tau_{E}^2\bigl(k y\bigl(2^{-n},
\theta_0,0\bigr)\bigr)\geq C_5 r_n^2.
\]
By the same proof of (\ref{Eq:UMC-low}) with some obvious
modifications, we have that
%
%
\begin{equation}
\label{ex2-14} \lim_{r\to0}\mathop{\sup_{y=y(r,\theta_0,0)}}_{
x,x+y\in I}
\frac{|X(x+y)-X(x)|} {r\sqrt{\log(1+r^{-1})}} \ge C >0.
\end{equation}
Reviewing the proof of Lemma~7.1.1 in Marcus and Rosen
\cite{MR-book}, one can easily get that
\[
\lim_{r\to0}\mathop{\sup_{y=y(r,\theta_0,0)}}_{ x,x+y\in I}
\frac{|X(x+y)-X(x)|} {r\sqrt{\log(1+r^{-1})}}\leq C, \qquad\mbox{a.s. for some constant } C<\infty
\]
implies that
\[
\lim_{r\to0}\mathop{\sup_{y=y(r,\theta_0,0)}}_{ x,x+y\in I}
\frac{|X(x+y)-X(x)|} {r\sqrt{\log(1+r^{-1})}}=C', \qquad\mbox{a.s. for some constant }
C'<\infty.
\]
Therefore, from Theorem~\ref{Th:UMC} and (\ref{ex2-14}) it follows
that
\[
\lim_{r\to0}\mathop{\sup_{y=y(r,\theta_0,0)}}_{
x,x+y\in I}
\frac{|X(x+y)-X(x)|} {r\sqrt{\log(1+r^{-1})}}=C\in(0, \infty).
\]
This and (\ref{ex2-7}) imply (\ref{ex2-15}).
\end{pf}

\section*{Acknowledgements}

Yuqiang Li's research is supported by Innovation Program
of Shanghai Municipal Education Commission (No:~13zz037), the
``Fundamental Research Funds for the Central Universities'' and
the 111 project (B14019).
Wensheng Wang's research is supported in part by
NSFC grant (No:~11071076).
Yimin Xiao's research is
supported in part by NSF grants DMS-1309856 and DMS-1307470.


%

\printhistory


\begin{thebibliography}{36}

\bibitem{A55}
%
\begin{barticle}[mr]
\bauthor{\bsnm{Anderson},~\bfnm{T.~W.}\binits{T.W.}}
(\byear{1955}).
\btitle{The integral of a symmetric unimodal function over a symmetric
convex set and some probability inequalities}.
\bjournal{Proc. Amer. Math. Soc.}
\bvolume{6}
\bpages{170--176}.
\bid{issn={0002-9939}, mr={0069229}}
\end{barticle}
%
\bptok{imsref}%
\endbibitem

\bibitem{AWX08}
%
\begin{barticle}[mr]
\bauthor{\bsnm{Ayache},~\bfnm{Antoine}\binits{A.}},
\bauthor{\bsnm{Wu},~\bfnm{Dongsheng}\binits{D.}} \AND
\bauthor{\bsnm{Xiao},~\bfnm{Yimin}\binits{Y.}}
(\byear{2008}).
\btitle{Joint continuity of the local times of fractional {B}rownian sheets}.
\bjournal{Ann. Inst. Henri Poincar\'e Probab. Stat.}
\bvolume{44}
\bpages{727--748}.
\bid{doi={10.1214/07-AIHP131}, issn={0246-0203}, mr={2446295}}
\end{barticle}
%
\bptok{imsref}%
\endbibitem

\bibitem{AX05}
%
\begin{barticle}[mr]
\bauthor{\bsnm{Ayache},~\bfnm{Antoine}\binits{A.}} \AND
\bauthor{\bsnm{Xiao},~\bfnm{Yimin}\binits{Y.}}
(\byear{2005}).
\btitle{Asymptotic properties and {H}ausdorff dimensions of fractional
{B}rownian sheets}.
\bjournal{J. Fourier Anal. Appl.}
\bvolume{11}
\bpages{407--439}.
\bid{doi={10.1007/s00041-005-4048-3}, issn={1069-5869}, mr={2169474}}
\end{barticle}
%
\bptok{imsref}%
\endbibitem

\bibitem{BM11}
%
\begin{barticle}[mr]
\bauthor{\bsnm{Baraka},~\bfnm{D.}\binits{D.}} \AND
\bauthor{\bsnm{Mountford},~\bfnm{T.~S.}\binits{T.S.}}
(\byear{2011}).
\btitle{The exact {H}ausdorff measure of the zero set of fractional
{B}rownian motion}.
\bjournal{J. Theoret. Probab.}
\bvolume{24}
\bpages{271--293}.
\bid{doi={10.1007/s10959-009-0271-1}, issn={0894-9840}, mr={2782718}}
\end{barticle}
%
\bptok{imsref}%
\endbibitem

\bibitem{BL02}
%
\begin{barticle}[mr]
\bauthor{\bsnm{Belinsky},~\bfnm{Eduard}\binits{E.}} \AND
\bauthor{\bsnm{Linde},~\bfnm{Werner}\binits{W.}}
(\byear{2002}).
\btitle{Small ball probabilities of fractional {B}rownian sheets via
fractional integration operators}.
\bjournal{J. Theoret. Probab.}
\bvolume{15}
\bpages{589--612}.
\bid{doi={10.1023/A:1016263614257}, issn={0894-9840}, mr={1922439}}
\end{barticle}
%
\bptok{imsref}%
\endbibitem

\bibitem{BMB06}
%
\begin{barticle}[auto:STB|2014/02/12|14:17:21]
\bauthor{\bsnm{Benson},~\bfnm{D.}\binits{D.}},
\bauthor{\bsnm{Meerschaert},~\bfnm{M.}\binits{M.}},
\bauthor{\bsnm{B{\"a}umer},~\bfnm{D.}\binits{D.}} \AND
\bauthor{\bsnm{Scheffler},~\bfnm{H.}\binits{H.}}
(\byear{2006}).
\btitle{Aquifer operator-scaling and the effect on solute mixing and
dispersion}.
\bjournal{Water Resour. Res.}
\bvolume{42}
\bpages{1--18}.
\end{barticle}
%
\bptok{imsref}%
\endbibitem

\bibitem{BL07}
%
\begin{barticle}[mr]
\bauthor{\bsnm{Bierm{\'e}},~\bfnm{Hermine}\binits{H.}} \AND
\bauthor{\bsnm{Lacaux},~\bfnm{C{\'e}line}\binits{C.}}
(\byear{2009}).
\btitle{H\"older regularity for operator scaling stable random fields}.
\bjournal{Stochastic Process. Appl.}
\bvolume{119}
\bpages{2222--2248}.
\bid{doi={10.1016/j.spa.2008.10.008}, issn={0304-4149}, mr={2531090}}
\end{barticle}
%
\bptok{imsref}%
\endbibitem

\bibitem{BMS07}
%
\begin{barticle}[mr]
\bauthor{\bsnm{Bierm{\'e}},~\bfnm{Hermine}\binits{H.}},
\bauthor{\bsnm{Meerschaert},~\bfnm{Mark~M.}\binits{M.M.}} \AND
\bauthor{\bsnm{Scheffler},~\bfnm{Hans-Peter}\binits{H.-P.}}
(\byear{2007}).
\btitle{Operator scaling stable random fields}.
\bjournal{Stochastic Process. Appl.}
\bvolume{117}
\bpages{312--332}.
\bid{doi={10.1016/j.spa.2006.07.004}, issn={0304-4149}, mr={2290879}}
\end{barticle}
%
\bptok{imsref}%
\endbibitem

\bibitem{BGT87}
%
\begin{bbook}[mr]
\bauthor{\bsnm{Bingham},~\bfnm{N.~H.}\binits{N.H.}},
\bauthor{\bsnm{Goldie},~\bfnm{C.~M.}\binits{C.M.}} \AND
\bauthor{\bsnm{Teugels},~\bfnm{J.~L.}\binits{J.L.}}
(\byear{1987}).
\btitle{Regular Variation}.
\bseries{Encyclopedia of Mathematics and Its Applications}
\bvolume{27}.
\blocation{Cambridge}:
\bpublisher{Cambridge Univ. Press}.
\bid{mr={0898871}}
\end{bbook}
%
\bptok{imsref}%
\endbibitem

\bibitem{BE03}
%
\begin{barticle}[mr]
\bauthor{\bsnm{Bonami},~\bfnm{Aline}\binits{A.}} \AND
\bauthor{\bsnm{Estrade},~\bfnm{Anne}\binits{A.}}
(\byear{2003}).
\btitle{Anisotropic analysis of some {G}aussian models}.
\bjournal{J. Fourier Anal. Appl.}
\bvolume{9}
\bpages{215--236}.
\bid{doi={10.1007/s00041-003-0012-2}, issn={1069-5869}, mr={1988750}}
\end{barticle}
%
\bptok{imsref}%
\endbibitem

\bibitem{CD99}
%
\begin{bbook}[mr]
\bauthor{\bsnm{Chil{\`e}s},~\bfnm{Jean-Paul}\binits{J.-P.}} \AND
\bauthor{\bsnm{Delfiner},~\bfnm{Pierre}\binits{P.}}
(\byear{1999}).
\btitle{Geostatistics: Modeling Spatial Uncertainty}.
\bseries{Wiley Series in Probability and Statistics: Applied
Probability and Statistics}.
\blocation{New York}:
\bpublisher{Wiley}.
\bid{doi={10.1002/9780470316993}, mr={1679557}}
\end{bbook}
%
\bptok{imsref}%
\endbibitem

\bibitem{DH99}
%
\begin{barticle}[mr]
\bauthor{\bsnm{Davies},~\bfnm{Steve}\binits{S.}} \AND
\bauthor{\bsnm{Hall},~\bfnm{Peter}\binits{P.}}
(\byear{1999}).
\btitle{Fractal analysis of surface roughness by using spatial data}.
\bjournal{J. R. Stat. Soc. Ser. B Stat. Methodol.}
\bvolume{61}
\bpages{3--37}.
\bnote{With discussion and a reply by the authors}.
\bid{doi={10.1111/1467-9868.00160}, issn={1369-7412}, mr={1664088}}
\bptnote{check related}%
\end{barticle}
%
\bptok{imsref}%
\endbibitem

\bibitem{Dunker2000}
%
\begin{barticle}[mr]
\bauthor{\bsnm{Dunker},~\bfnm{Thomas}\binits{T.}}
(\byear{2000}).
\btitle{Estimates for the small ball probabilities of the fractional
{B}rownian sheet}.
\bjournal{J. Theoret. Probab.}
\bvolume{13}
\bpages{357--382}.
\bid{doi={10.1023/A:1007897525164}, issn={0894-9840}, mr={1777539}}
\end{barticle}
%
\bptok{imsref}%
\endbibitem

\bibitem{K96}
%
\begin{barticle}[mr]
\bauthor{\bsnm{Kamont},~\bfnm{Anna}\binits{A.}}
(\byear{1996}).
\btitle{On the fractional anisotropic {W}iener field}.
\bjournal{Probab. Math. Statist.}
\bvolume{16}
\bpages{85--98}.
\bid{issn={0208-4147}, mr={1407935}}
\end{barticle}
%
\bptok{imsref}%
\endbibitem

\bibitem{KL2002}
%
\begin{barticle}[mr]
\bauthor{\bsnm{K{\"u}hn},~\bfnm{Thomas}\binits{T.}} \AND
\bauthor{\bsnm{Linde},~\bfnm{Werner}\binits{W.}}
(\byear{2002}).
\btitle{Optimal series representation of fractional {B}rownian sheets}.
\bjournal{Bernoulli}
\bvolume{8}
\bpages{669--696}.
\bid{issn={1350-7265}, mr={1935652}}
\end{barticle}
%
\bptok{imsref}%
\endbibitem

\bibitem{LX09}
%
\begin{barticle}[mr]
\bauthor{\bsnm{Li},~\bfnm{Yuqiang}\binits{Y.}} \AND
\bauthor{\bsnm{Xiao},~\bfnm{Yimin}\binits{Y.}}
(\byear{2011}).
\btitle{Multivariate operator-self-similar random fields}.
\bjournal{Stochastic Process. Appl.}
\bvolume{121}
\bpages{1178--1200}.
\bid{doi={10.1016/j.spa.2011.02.005}, issn={0304-4149}, mr={2794972}}
\end{barticle}
%
\bptok{imsref}%
\endbibitem

\bibitem{LX13}
%
\begin{barticle}[mr]
\bauthor{\bsnm{Li},~\bfnm{Yuqiang}\binits{Y.}} \AND
\bauthor{\bsnm{Xiao},~\bfnm{Yimin}\binits{Y.}}
(\byear{2013}).
\btitle{A class of fractional {B}rownian fields from branching systems
and their regularity properties}.
\bjournal{Infin. Dimens. Anal. Quantum Probab. Relat. Top.}
\bvolume{16}
\bpages{1350023, 33}.
\bid{doi={10.1142/S0219025713500239}, issn={0219-0257}, mr={3125843}}
\end{barticle}
%
\bptok{imsref}%
\endbibitem

\bibitem{LuanXiao10}
%
\begin{barticle}[mr]
\bauthor{\bsnm{Luan},~\bfnm{Nana}\binits{N.}} \AND
\bauthor{\bsnm{Xiao},~\bfnm{Yimin}\binits{Y.}}
(\byear{2012}).
\btitle{Spectral conditions for strong local nondeterminism and exact
{H}ausdorff measure of ranges of {G}aussian random fields}.
\bjournal{J. Fourier Anal. Appl.}
\bvolume{18}
\bpages{118--145}.
\bid{doi={10.1007/s00041-011-9193-2}, issn={1069-5869}, mr={2885561}}
\end{barticle}
%
\bptok{imsref}%
\endbibitem

\bibitem{MR-book}
%
\begin{bbook}[mr]
\bauthor{\bsnm{Marcus},~\bfnm{Michael~B.}\binits{M.B.}} \AND
\bauthor{\bsnm{Rosen},~\bfnm{Jay}\binits{J.}}
(\byear{2006}).
\btitle{Markov Processes, {G}aussian Processes, and Local Times}.
\bseries{Cambridge Studies in Advanced Mathematics}
\bvolume{100}.
\blocation{Cambridge}:
\bpublisher{Cambridge Univ. Press}.
\bid{doi={10.1017/CBO9780511617997}, mr={2250510}}
\end{bbook}
%
\bptok{imsref}%
\endbibitem

\bibitem{MS2001}
%
\begin{barticle}[mr]
\bauthor{\bsnm{Mason},~\bfnm{David~M.}\binits{D.M.}} \AND
\bauthor{\bsnm{Shi},~\bfnm{Zhan}\binits{Z.}}
(\byear{2001}).
\btitle{Small deviations for some multi-parameter {G}aussian processes}.
\bjournal{J.~Theoret. Probab.}
\bvolume{14}
\bpages{213--239}.
\bid{doi={10.1023/A:1007833401562}, issn={0894-9840}, mr={1822902}}
\end{barticle}
%
\bptok{imsref}%
\endbibitem

\bibitem{RVbook}
%
\begin{bbook}[mr]
\bauthor{\bsnm{Meerschaert},~\bfnm{Mark~M.}\binits{M.M.}} \AND
\bauthor{\bsnm{Scheffler},~\bfnm{Hans-Peter}\binits{H.-P.}}
(\byear{2001}).
\btitle{Limit Distributions for Sums of Independent Random Vectors: Heavy Tails in Theory and Practice}.
\blocation{New York}:
\bpublisher{Wiley}.
\bid{mr={1840531}}
\end{bbook}
%
\bptok{imsref}%
\endbibitem

\bibitem{MWX10}
%
\begin{barticle}[mr]
\bauthor{\bsnm{Meerschaert},~\bfnm{Mark~M.}\binits{M.M.}},
\bauthor{\bsnm{Wang},~\bfnm{Wensheng}\binits{W.}} \AND
\bauthor{\bsnm{Xiao},~\bfnm{Yimin}\binits{Y.}}
(\byear{2013}).
\btitle{Fernique-type inequalities and moduli of continuity for
anisotropic {G}aussian random fields}.
\bjournal{Trans. Amer. Math. Soc.}
\bvolume{365}
\bpages{1081--1107}.
\bid{doi={10.1090/S0002-9947-2012-05678-9}, issn={0002-9947}, mr={2995384}}
\end{barticle}
%
\bptok{imsref}%
\endbibitem

\bibitem{ST94}
%
\begin{bbook}[mr]
\bauthor{\bsnm{Samorodnitsky},~\bfnm{Gennady}\binits{G.}} \AND
\bauthor{\bsnm{Taqqu},~\bfnm{Murad~S.}\binits{M.S.}}
(\byear{1994}).
\btitle{Stable Non-{G}aussian Random Processes: Stochastic Models with Infinite Variance}.
\bseries{Stochastic Modeling}.
\blocation{New York}:
\bpublisher{Chapman \& Hall}.
\bid{mr={1280932}}
\end{bbook}
%
\bptok{imsref}%
\endbibitem

\bibitem{Stein05}
%
\begin{barticle}[mr]
\bauthor{\bsnm{Stein},~\bfnm{Michael~L.}\binits{M.L.}}
(\byear{2005}).
\btitle{Space-time covariance functions}.
\bjournal{J. Amer. Statist. Assoc.}
\bvolume{100}
\bpages{310--321}.
\bid{doi={10.1198/016214504000000854}, issn={0162-1459}, mr={2156840}}
\end{barticle}
%
\bptok{imsref}%
\endbibitem

\bibitem{Stein12}
%
\begin{barticle}[mr]
\bauthor{\bsnm{Stein},~\bfnm{Michael~L.}\binits{M.L.}}
(\byear{2013}).
\btitle{On a class of space-time intrinsic random functions}.
\bjournal{Bernoulli}
\bvolume{19}
\bpages{387--408}.
\bid{doi={10.3150/11-BEJ405}, issn={1350-7265}, mr={3037158}}
\end{barticle}
%
\bptok{imsref}%
\endbibitem

\bibitem{T95}
%
\begin{barticle}[mr]
\bauthor{\bsnm{Talagrand},~\bfnm{Michel}\binits{M.}}
(\byear{1995}).
\btitle{Hausdorff measure of trajectories of multiparameter fractional
{B}rownian motion}.
\bjournal{Ann. Probab.}
\bvolume{23}
\bpages{767--775}.
\bid{issn={0091-1798}, mr={1334170}}
\end{barticle}
%
\bptok{imsref}%
\endbibitem

\bibitem{Wack98}
%
\begin{bbook}[auto:STB|2014/02/12|14:17:21]
\bauthor{\bsnm{Wackernagel},~\bfnm{H.}\binits{H.}}
(\byear{1998}).
\btitle{Multivariate Geostatistics: An Introduction with Applications}.
\blocation{New York}:
\bpublisher{Springer}.
\end{bbook}
%
\bptok{imsref}%
\endbibitem

\bibitem{Wang07}
%
\begin{barticle}[mr]
\bauthor{\bsnm{Wang},~\bfnm{Wensheng}\binits{W.}}
(\byear{2007}).
\btitle{Almost-sure path properties of fractional {B}rownian sheet}.
\bjournal{Ann. Inst. Henri Poincar\'e Probab. Stat.}
\bvolume{43}
\bpages{619--631}.
\bid{doi={10.1016/j.anihpb.2006.09.005}, issn={0246-0203}, mr={2347099}}
\end{barticle}
%
\bptok{imsref}%
\endbibitem

\bibitem{WX07}
%
\begin{barticle}[mr]
\bauthor{\bsnm{Wu},~\bfnm{Dongsheng}\binits{D.}} \AND
\bauthor{\bsnm{Xiao},~\bfnm{Yimin}\binits{Y.}}
(\byear{2007}).
\btitle{Geometric properties of fractional {B}rownian sheets}.
\bjournal{J. Fourier Anal. Appl.}
\bvolume{13}
\bpages{1--37}.
\bid{doi={10.1007/s00041-005-5078-y}, issn={1069-5869}, mr={2296726}}
\end{barticle}
%
\bptok{imsref}%
\endbibitem

\bibitem{WX11}
%
\begin{barticle}[mr]
\bauthor{\bsnm{Wu},~\bfnm{Dongsheng}\binits{D.}} \AND
\bauthor{\bsnm{Xiao},~\bfnm{Yimin}\binits{Y.}}
(\byear{2011}).
\btitle{On local times of anisotropic {G}aussian random fields}.
\bjournal{Commun. Stoch. Anal.}
\bvolume{5}
\bpages{15--39}.
\bid{issn={0973-9599}, mr={2808534}}
\end{barticle}
%
\bptok{imsref}%
\endbibitem

\bibitem{X97a}
%
\begin{barticle}[mr]
\bauthor{\bsnm{Xiao},~\bfnm{Yimin}\binits{Y.}}
(\byear{1997}).
\btitle{Hausdorff measure of the graph of fractional {B}rownian motion}.
\bjournal{Math. Proc. Cambridge Philos. Soc.}
\bvolume{122}
\bpages{565--576}.
\bid{doi={10.1017/S0305004197001783}, issn={0305-0041}, mr={1466658}}
\end{barticle}
%
\bptok{imsref}%
\endbibitem

\bibitem{X97b}
%
\begin{barticle}[mr]
\bauthor{\bsnm{Xiao},~\bfnm{Yimin}\binits{Y.}}
(\byear{1997}).
\btitle{H\"older conditions for the local times and the {H}ausdorff
measure of the level sets of {G}aussian random fields}.
\bjournal{Probab. Theory Related Fields}
\bvolume{109}
\bpages{129--157}.
\bid{doi={10.1007/s004400050128}, issn={0178-8051}, mr={1469923}}
\end{barticle}
%
\bptok{imsref}%
\endbibitem

\bibitem{Xiao09}
%
\begin{bincollection}[mr]
\bauthor{\bsnm{Xiao},~\bfnm{Yimin}\binits{Y.}}
(\byear{2009}).
\btitle{Sample path properties of anisotropic {G}aussian random fields}.
In \bbooktitle{A Minicourse on Stochastic Partial Differential Equations}.
\bseries{Lecture Notes in Math.}
\bvolume{1962}
\bpages{145--212}.
\bpublisher{Springer, Berlin}.
\bid{doi={10.1007/978-3-540-85994-9_5}, mr={2508776}}
\end{bincollection}
%
\bptok{imsref}%
\endbibitem

\bibitem{Xiao10}
%
\begin{barticle}[mr]
\bauthor{\bsnm{Xiao},~\bfnm{Yimin}\binits{Y.}}
(\byear{2010}).
\btitle{Uniform modulus of continuity of random fields}.
\bjournal{Monatsh. Math.}
\bvolume{159}
\bpages{163--184}.
\bid{doi={10.1007/s00605-009-0133-z}, issn={0026-9255}, mr={2564392}}
\end{barticle}
%
\bptok{imsref}%
\endbibitem

\bibitem{XZ04}
%
\begin{barticle}[mr]
\bauthor{\bsnm{Xiao},~\bfnm{Yimin}\binits{Y.}} \AND
\bauthor{\bsnm{Zhang},~\bfnm{Tusheng}\binits{T.}}
(\byear{2002}).
\btitle{Local times of fractional {B}rownian sheets}.
\bjournal{Probab. Theory Related Fields}
\bvolume{124}
\bpages{204--226}.
\bid{doi={10.1007/s004400200210}, issn={0178-8051}, mr={1936017}}
\end{barticle}
%
\bptok{imsref}%
\endbibitem

\bibitem{Zhang07}
%
\begin{barticle}[mr]
\bauthor{\bsnm{Zhang},~\bfnm{Hao}\binits{H.}}
(\byear{2007}).
\btitle{Maximum-likelihood estimation for multivariate spatial linear
coregionalization models}.
\bjournal{Environmetrics}
\bvolume{18}
\bpages{125--139}.
\bid{doi={10.1002/env.807}, issn={1180-4009}, mr={2345650}}
\end{barticle}
%
\bptok{imsref}%
\endbibitem

\end{thebibliography}
\end{document}